\documentclass[pamm,a4paper,fleqn]{w-art}
\usepackage{times,cite,w-thm}
\usepackage[T1]{fontenc}
\usepackage[utf8]{inputenc}
\usepackage{graphicx}
\usepackage{siunitx}
\usepackage{placeins}
\usepackage{amsmath,amsthm,amsfonts,amssymb}
\usepackage{mathtools}
\usepackage{booktabs}
\usepackage{algorithm}
\usepackage[noend]{algorithmic}
\renewcommand{\algorithmiccomment}[1]{\bgroup\hfill$\blacktriangleright$~#1\egroup}
\usepackage{hyperref}
\newcommand{\doi}[1]{\href{https://doi.org/#1}{doi:#1}}

\newcommand{\arXiv}[1]{\href{https://arxiv.org/abs/#1}{arXiv:#1}}
\newtheorem{Problem}[theorem]{Problem}
\usepackage{tikz}
\usepackage{tikz-3dplot}
\usetikzlibrary{matrix,calc}
\usepackage{pgf,pgfplots}
\usepackage{subcaption}
\usepgfplotslibrary{groupplots,colorbrewer}
\usepgfplotslibrary{colormaps} \pgfplotsset{ colormap/Set1-5, cycle multiindex*
  list={ mark list*\nextlist Set1-5\nextlist }, every axis/.append style =
  {thick},%
}%
\pgfplotsset{tick style = {thick,black}}%
\newcommand{\functionforplot}[2]{(-#1^2 + 2*#1 + #2^2 + #1*#2/4 -#2 + 2)}
\newif\ifMAKEPICS
\MAKEPICSfalse

\ifMAKEPICS
\usepackage[cleanup,subfolder]{gnuplottex}
\usepackage{xparse}

\ExplSyntaxOn
\DeclareExpandableDocumentCommand{\convertlen}{ O{cm} m }
{
  \dim_to_decimal_in_unit:nn { #2 } { 1 #1 } cm
}
\ExplSyntaxOff
\fi
\begin{document}

\TitleLanguage[EN]
\title[The short title]{Multi-goal-oriented anisotropic error control and mesh adaptivity for time-dependent convection-dominated problems}

\author{\firstname{Markus} \lastname{Bause}\inst{1,}%
\footnote{e-mail \ElectronicMail{bause@hsu-hh.de}}}
\address[\inst{1}]{\CountryCode[DE]Helmut Schmidt University, University of the German Federal Armed Forces Hamburg, Faculty of mechanical and Civil Engineering, Chair of Numerical Mathematics, Holstenhofweg 85, 22043 Hamburg, Germany}
\author{\firstname{Marius Paul} \lastname{Bruchh\"auser}\inst{1,}%
\footnote{Corresponding author: e-mail \ElectronicMail{bruchhaeuser@hsu-hh.de}, phone +49\,40\,6541\,2739, fax +49\,40\,6541\,2690}}
\author{\firstname{Bernhard} \lastname{Endtmayer}\inst{2,3,}%
     \footnote{e-mail \ElectronicMail{endtmayer@ifam.uni-hannover.de}}}
\address[\inst{2}]{\CountryCode[DE]Leibniz Universität Hannover, Institut für Angewandte Mathematik, Welfengarten 1, 30167 Hannover, Germany}
\address[\inst{3}]{\CountryCode[DE]Cluster of Excellence PhoenixD (Photonics, Optics, and Engineering – Innovation Across
Disciplines), Leibniz Universität Hannover, Germany}
\author{\firstname{Nils} \lastname{Margenberg}\inst{1,}%
     \footnote{e-mail \ElectronicMail{margenbn@hsu-hh.de}}}
\author{\firstname{Ioannis} \lastname{Toulopoulos}\inst{4,}%
     \footnote{itoulopoulos@uowm.gr}}
\address[\inst{4}]{\CountryCode[EN]University of Western Macedonia, Department of Informatics, Fourka Area, 52100 Kastoria, Greece}
\author{\firstname{Thomas} \lastname{Wick}\inst{2,3,}%
     \footnote{e-mail \ElectronicMail{wick@ifam.uni-hannover.de}}}
\AbstractLanguage[EN]
\begin{abstract}
In this work, we present an anisotropic multi-goal error control based on the Dual Weighted Residual (DWR) method for time-dependent convection-diffusion-reaction (CDR) equations.
This multi-goal oriented approach allows for an accurate and efficient error control with regard to several quantities of interest simultaneously.
Using anisotropic interpolation and restriction operators, we obtain  elementwise error indicators in space and time, where the spatial indicators are additionally separated with respect to the single directions.
The directional error indicators quantify anisotropy of the solution with
respect to the goals, and produce adaptive, anisotropic meshes that efficiently capture layers.
To prevent spurious oscillations the streamline upwind Petrov-Galerkin (SUPG) method is applied to stabilize the underlying system in the case of high P\'{e}clet numbers.
Numerical examples show efficiency and robustness of the proposed approach for
several goal quantities using established benchmarks for convection-dominated transport.
\end{abstract}
\maketitle                   

\section{Introduction}
\label{sec:1:intro}
In the recent past, adaptive methods have become an important additional tool to stabilization in order to handle the challenges associated with an accurate and efficient nunmerical solution for convection-dominated problems, cf., e.g., \cite{BruchhaeuserCCY24,BruchhaeuserAJK23,BruchhaeuserBB24,BruchhaeuserBSB19,BruchhaeuserABR17,BruchhaeuserYK13,BruchhaeuserBK13}.
The solutions of such problems are typically characterized by sharp moving fronts and interior or boundary layers.
Since the pioneering works of the 1980's (cf., e.g.,~\cite{BruchhaeuserBH82,BruchhaeuserHMM86}) a lot of effort was spent on  the development of appropriate stabilization techniques to avoid spurious, unphysical oscillations.
This includes both residual-based stabilization as well as algebraic stabilization methods. For a small review of the most common stabilization techniques, we refer to, e.g., \cite{BruchhaeuserRST08,BruchhaeuserJKN18}.
As shown in a comparative study for time-dependent convection-diffusion-reaction equations in \cite{BruchhaeuserJS08}, conventional stabilization techniques based on standard meshes fail to avoid oscillations or exhibit so-called smearing effects even after a careful fine-tuning of the respective stabilization parameters.
Similar results have been perceived in three dimensions; cf. \cite{BruchhaeuserJS09}.

Recently, the combination of stabilization and adaptivity applied to CDR equations was extended by means of anisotropic mesh refinement, cf.~\cite{BruchhaeuserKS25,BruchhaeuserBBEMTW25}.
Anisotropic refinement along dominant error directions may improve computational efficiency significantly. This approach has been successfully applied to many fields of partial differential equations in general, cf., e.g., \cite{BruchhaeuserZBHKD22,BruchhaeuserCHDG20,BruchhaeuserR13,BruchhaeuserLH10,BruchhaeuserP03} and in particular to CDR equations, cf., e.g., \cite{BruchhaeuserYD12,BruchhaeuserGHH08,BruchhaeuserFPZ01,BruchhaeuserAL96,BruchhaeuserKR90}.
For a general overview of anisotropic finite elements we refer to the
monograph of Apel \cite{BruchhaeuserA99}.

In this work, we present an anisotropic multi-goal error control based on the Dual Weighted Residual method \cite{BruchhaeuserBR98,BruchhaeuserBR01,BruchhaeuserBR03} for time-dependent convection-diffusion-reaction equations.
Motivated by ideas of Richter~\cite{BruchhaeuserR10,BruchhaeuserR13}, we combine goal-oriented error estimation with anisotropic refinement.
Our proposed method provides directional error indicators in space and time with respect to user-defined goals. Thereby, the directional separation of elementwise indicators naturally leads to anisotropic mesh refinement.
Often in numerical simulations, not the entire numerical solution is of interest, but certain quantities of interest.
This multi-goal oriented approach allows for an accurate and efficient error control with regard to several quantities of interest simultaneously.
In 2003, goal-oriented a posteriori error estimation was first proposed for multiple quantities of interest in \cite{BruchhaeuserHH03}.
Since then, multi-goal-oriented error control has been successfully applied to several mathematical model problems including inter alia elliptic, nonlinear, aerodynamic flow, coupled problems and many more, cf., e.g.,
\cite{BruchhaeuserH08,BruchhaeuserEW17,beuchler2024mathematical,BruchhaeuserP10}.
For a briefly review to multi-goal error control we refer to, e.g., \cite{BruchhaeuserELRSW24}.

This work is organized as follows.
In Sect.\ref{sec:2:model-disc} we introduce our model problem as well as the stabilized space-time discretizations.
Sec.~\ref{sec:error_rep} derives an a posteriori error representation with respect to a single goal.
The main focus is on Sect.\ref{sec:anisotropic-method}, where we introduce our anisotropic multi-goal oriented error control approach, present the underlying space-time adaptive algorithm and give insight into some practical aspects of the implementation.
Finally, Sec.~\ref{sec:5:numerical_examples} illustrates and validates our approach by means of typical numerical benchmarks for convection-dominated problems.

\section{Model Problem and Discretization in Space and Time}
\label{sec:2:model-disc}
In this work, we study the following time-dependent convection-diffusion-reaction equation:
\begin{equation}
\label{eq:1:CDRoriginal}
\begin{array}{r@{\;}c@{\;}l@{\;} @{\,\,}l @{\,\,}l @{\,}l}
\partial_t u
- \nabla \cdot  (\varepsilon \nabla u)
+ \boldsymbol{b} \cdot \nabla u
+ \alpha u & = & f & \mbox{in } & \mathcal{Q} & = \Omega \times I \,,
\\
 u & = & u_D & \mbox{on } & \Sigma_D & = \Gamma_D \times I\,,
\\
 \varepsilon \nabla u \cdot \boldsymbol{n} & = & u_N & \mbox{on } & \Sigma_N & = \Gamma_N \times I\,,
\\
u(0) & = & u_{0} & \mbox{on } & \Sigma_0 & = \Omega\times \{0\} \,,
\end{array}
\end{equation}
in the space-time domain $\mathcal{Q}$, where $\Omega\subset\mathbb{R}^d$,
with $d=2$ or $d=3$, is a polygonal or polyhedral bounded domain with Lipschitz
boundary $\partial\Omega$ and $I=(0,T), 0 < T < \infty$, is a finite
time interval.
Here, $\partial\Omega = \Gamma_D \cup \Gamma_N\,,\Gamma_D \neq
\emptyset$ denotes the partition of the boundary with outer unit normal vector
$\boldsymbol{n}$, where $\Gamma_D$ denotes the
Dirichlet part and $\Gamma_N$ the Neumann part, respectively.
Furthermore, let $V:=\big\{v \in H^1(\Omega)|v_{|\Gamma_D}=0\big\}$ and $V'$ denotes the adjoint space of $V$.
To ensure the well-posedness of Eq.~\eqref{eq:1:CDRoriginal} we
assume that $0 < \varepsilon \leq 1$ is a constant diffusion coefficient,
$\boldsymbol{b} \in L^{\infty}\big(I;W^{1,\infty}(\Omega)^d\big)$
is the flow field or convection field,
$\alpha \in L^{\infty}\big(I;L^{\infty}(\Omega)\big)$
is a non-negative ($\alpha \geq 0$) reaction coefficient,
$u_0\in L^2(\Omega)$
is a given initial condition,
$f \in L^{2}(I;V')$ is a given source of the unknown scalar quantity $u$,
$u_D \in L^{2}(I;H^{\frac{1}{2}}(\Gamma_D))$ is a given function specifying the
Dirichlet boundary condition, and $u_N \in L^{2}(I;H^{-\frac{1}{2}}(\Gamma_N))$ is a given
function specifying the Neumann boundary condition.
Furthermore, it will be assumed that either $\nabla \cdot \boldsymbol{b}(\boldsymbol{x},t) = 0$
and $\alpha(\boldsymbol{x},t) \geq 0$, or there exists a
positive constant $c_0$ such that
$
\alpha(\boldsymbol{x},t)-\frac{1}{2} \textnormal{div}\;\boldsymbol{b}(\boldsymbol{x},t) \geq c_0 > 0
\;\;\forall (\boldsymbol{x},t) \in \bar{\Omega}\times \bar{I}\,,
$
which are standard assumptions for convection-dominated equations of type~\eqref{eq:1:CDRoriginal}, cf., e.g.,~\cite{BruchhaeuserAJ15,BruchhaeuserRST08}.

Henceforth, for the sake of simplicity, we deal with homogeneous Dirichlet
boundary values $u_D = 0$ on $\Gamma_D=\partial\Omega$ only. This implies that here $V\coloneq H_0^1(\Omega)$. In the numerical examples in Sec.~\ref{sec:5:numerical_examples}, we also
consider more general boundary conditions, which are incorporated as described in~\cite[Ch.~3.3]{BruchhaeuserBR03}.
It is well known that problem \eqref{eq:1:CDRoriginal} along with the above
conditions admits a unique weak solution $u \in \mathcal{V}\coloneqq \big\{v\in L^{2}\big(I; V\big)\big|\;\partial_{t}v\in L^{2}(I;V')\big\}\,,$ that satisfies the following variational formulation; cf., e.g.~\cite{BruchhaeuserRST08,BruchhaeuserJKN18}.
\begin{Problem}
Find $u \in \mathcal{V}$ such that
\begin{equation}
\label{eq:2:WeakCDRsteady}
A(u)(\varphi) = F(\varphi) \quad \forall \varphi \in \mathcal{V}\,,
\end{equation}
where the bilinear form $A:\mathcal{V}\times \mathcal{V} \rightarrow \mathbb{R}$ and the linear form $F:\mathcal{V}\rightarrow \mathbb{R}$ are defined by
\begin{eqnarray}
\label{eq:3:BilinearformA}
A(u)(\varphi) & \coloneqq & \int_{I}\big\{(\partial_t u,\varphi)
+a(u)(\varphi)\big\} \mathrm{d} t + (u(0),\varphi(0))\,,
\\[1ex]
\label{eq:4:LinearformF}
F(\varphi) & \coloneqq & \int_I(f,\varphi)\;\mathrm{d}t + (u_0,\varphi(0))\,,
\end{eqnarray}
with the inner bilinear form $a:V \times V \rightarrow \mathbb{R}$, given by
\begin{equation}
\label{eq:5:aBilinearform}
a(u)(\varphi)\coloneqq (\varepsilon\nabla u, \nabla \varphi)
+(\boldsymbol{b}\cdot \nabla u, \varphi) + (\alpha u,\varphi)\,.
\end{equation}
\end{Problem}
We note that the initial condition is incorporated into the variational problem. The weak formulation, given by Eq.~\eqref{eq:2:WeakCDRsteady}, is now the starting point for the variational discretization in space and time using Galerkin finite element methods.
\subsection{Discretization in Time}
\label{sec:2.1:disc_time}
For the discretization in time we use a discontinuous Galerkin method $dG(r)$
with an arbitrary polynomial degree $r\geq0$.
Let $\mathcal T_{\tau}$ be a partition of the closure of the time domain
$\bar{I}=[0,T]$ into left-open subintervals $I_n\coloneqq (t_{n-1},t_n]$, $n=1,\dots,N$,
with $0=:t_0<t_1<\dots<t_N\coloneqq T$ and time step sizes $\tau_n=t_n-t_{n-1}$ and the global time discretization
parameter $\tau=\max_{n}\,\tau_{n}$.
Since the set of time intervals \( I_n \) is finite, it is natural to partition the global space-time cylinder \( \mathcal{Q} = \Omega \times I \) into space-time slabs defined as \( \mathcal{Q}_n = \Omega \times I_n \).
On the subintervals $I_n$, we define the time-discrete function space
$\mathcal{V}_{\tau}^{r}$
\begin{equation}
\label{eq:6:Def_V_tau_dGr}
 \begin{aligned}
\mathcal{V}_{\tau}^{r} \coloneqq
 \Big\{u_{\tau}\in L^{2}(I; V)\big|
 u_{\tau}|_{I_{n}}\in \mathcal{P}_{r}(I_{n}; V),
 u_{\tau}(0)\in L^2(\Omega), n=1,\dots,N
\Big\}\,,
\end{aligned}
\end{equation}
where $\mathcal{P}_{r}(\bar{I}_{n}; V)$ denotes the space of all
polynomials in time up to degree $r\geq0$ on $I_n$ with values in $V\,.$
For some discontinuous in time function $u_{\tau}\in \mathcal{V}_{\tau}^{r}$
we define the limits $u_{\tau,n}^{\pm}$ from above and below of $u_{\tau}$ at
$t_n$ as well as their jump at $t_n$ by
\begin{displaymath}
\begin{array}{lcrclcr}
u_{\tau,n}^{\pm}
& \coloneqq &
\displaystyle\lim_{t\to t_n\pm0} u_\tau(t) \,,
&
[u_{\tau}]_{n} & \coloneqq & u_{\tau,n}^{+}
-u_{\tau,n}^{-} \,.
\end{array}
\end{displaymath}
Then, the semi-discrete in time scheme of Eq.~\eqref{eq:2:WeakCDRsteady} reads
as follows:
\textit{
Find $u_{\tau} \in \mathcal{V}_{\tau}^{r}$ such that
}
\begin{equation}
\label{eq:7:dGDiscTime}
A_{\tau}(u_{\tau})(\varphi_{\tau}) =
F_{\tau}(\varphi_{\tau}) \quad \forall \varphi_{\tau}
\in \mathcal{V}_{\tau}^{r}\,,
\end{equation}
\textit{
where the bilinear form $A_{\tau}(\cdot)(\cdot)$ and the linear form
$F_{\tau}(\cdot)$ are defined by
}
\begin{eqnarray}
\label{eq:8:BilinearFormAtau}
A_{\tau}(u_{\tau})(\varphi_{\tau}) & \coloneqq &
\displaystyle\sum_{n=1}^{N}\int_{I_n}
\big\{(\partial_t u_{\tau},\varphi_{\tau})
+a(u_{\tau})(\varphi_{\tau})
\big\} \mathrm{d} t
+ \displaystyle\sum_{n=2}^N
([u_{\tau}]_{n-1},\varphi_{\tau,n-1}^+ )
+ (u_{\tau,0}^{+},\varphi_{\tau,0}^{+})\,,
\\[1ex]
\label{eq:9:LinearFormFtau}
F_{\tau}(\varphi_{\tau}) & \coloneqq & \displaystyle\int_I(f,\varphi)\;\mathrm{d}t
+(u_0,\varphi_{\tau,0}^{+})\,,
\end{eqnarray}
with the inner bilinear form $a(\cdot)(\cdot)$ being defined by
Eq.~\eqref{eq:5:aBilinearform}.
\subsection{Discretization in Space and Stabilization}
\label{sec:3.3:disc_space}
Next, we describe the Galerkin finite element approximation in space for the
semi-discrete time scheme~\eqref{eq:7:dGDiscTime}. We use Lagrange-type finite
element spaces of continuous, piecewise polynomial functions. The spatial
discretization is based on a decomposition \( \mathcal{T}_h \) of the domain
\( \Omega \) into disjoint elements \( K \), such that
\( \overline{\Omega} = \cup_{K\in\mathcal{T}_h} \overline{K} \). For
\( d=2,3 \), we use quadrilateral and hexahedral elements, respectively. Each
element \( K \in \mathcal{T}_h \) is mapped from the reference element
$\hat{K}=(0,1)^d$ via an iso-parametric transformation \( \boldsymbol T_K: \hat{K} \to K \) satisfying \( \det(\boldsymbol T_K)(\hat{x}) > 0 \) for all \( \hat{x} \in (0,1)^d \).
Following~\cite{BruchhaeuserR10}, we decompose
\begin{equation}
\label{eq:TK}
\boldsymbol T_K \coloneq \boldsymbol R_K \circ\boldsymbol S_{c,K} \circ\boldsymbol S_{h,K} \circ\boldsymbol P_K,
\end{equation}
where \(\boldsymbol R_K \) is a rotation and translation, \(\boldsymbol S_{c,K} \) an anisotropic scaling, \(\boldsymbol S_{h,K} \) a shearing, and \(\boldsymbol P_K \) is a nonlinear component.
To account for anisotropic elements, we relax the standard shape-regularity
assumptions and require only uniform boundedness of \(\boldsymbol S_{h,K} \) and \(\boldsymbol P_K \)
for all \( K \in \mathcal{T}_h \). The element diameter is denoted by \( h_K \),
with the global discretization parameter defined as
\( h \coloneqq \max_{K\in\mathcal{T}_h} h_K \).

Our mesh adaptation procedure employs local refinement with hanging nodes.
Global conformity of the finite element approach is preserved by eliminating
degrees of freedom at hanging nodes via interpolation between neighboring
regular nodes; see~\cite[Ch.~4.2]{BruchhaeuserBR03} and~\cite{BruchhaeuserCO84}.
\begin{definition}
On a subset $\mathcal Z_h\subseteq \mathcal T_h$, we define the discrete finite element space
$$
V_h^{p}(\mathcal Z_h)\coloneqq
\big\{v\in C(\overline{\Omega})\mid v_{|K}
\in Q_{p}(K)\,,\forall K\in\mathcal Z_h
\big\}\cap V\,,
$$
where, $Q_{p}(K)$ is the mapped finite element from
\[
\hat{Q}_{p}(\hat{K}) \coloneqq \bigotimes_{i=1}^{d} \mathcal{P}_{p}([0,1])\,,
\]
by~\eqref{eq:TK}, where \(\mathcal{P}_{p}([0,1])\) is the space of univariate
polynomials of degree $p$. On $\mathcal{T}_h$ we put
$$V_h^{p}=V_h^{p}(\mathcal{T}_h)\,.$$
\end{definition}
In order to split the error in directional contributions, we further
introduce finite elements with anisotropic polynomial degree.
\begin{definition}
\label{def:anisoFEM}
Let \((p_1,\dots,p_d)\in\mathbb{N}^d\) be a multi-index.
On the reference element $\hat{K}$ the anisotropic polynomial space is defined as
\[
\hat{Q}_{p_1,\dots,p_d}(\hat{K}) \coloneqq \bigotimes_{i=1}^{d} \mathcal{P}_{p_i}([0,1]),
\]
where \(\mathcal{P}_{p_i}([0,1])\) denotes the space of univariate polynomials
of degree at most \(p_i\).
\end{definition}
In Def.~\ref{def:anisoFEM}, the DoFs are associated with Gauss--Lobatto nodes.
The total number of DoFs per element is
\(
N_K^{p_1,\dots,p_d} = \prod_{i=1}^{d} (p_i+1).
\)
We note that this is a straightforward generalization of the isotropic case, where for
\(p_1=\cdots=p_d \coloneqq p\), we put $\hat{Q}_{p}$.
For \(i\in\{1,\dots,d\}\) and given integers \(p\) and \(q\), we define
\[
\hat{Q}^{p,q}_i \coloneqq \hat{Q}^{p_1,\dots,p_d}, \quad
p_j =
\begin{cases}
  p, & j\neq i \\
  q, & j=i
\end{cases}\,.
\]
Further, we denote by $N_i^{p,q}$ the set of DoFs on the reference element
$\hat K$.

Let the fully discrete function space be given by
\begin{equation}
\label{eq:10:Def_V_tau_h_dGr_p}
\mathcal{V}_{\tau h}^{r,p} \coloneqq \Big\{
u_{\tau h}\in X_{\tau}^{r} \big|
u_{\tau h|{I_n}} \in \mathcal{P}_r(I_n;V_h^{p})
\,, u_{\tau h}(0) \in V_h^{p},
n=1,\dots,N
\Big\}
\subseteq L^{2}(I; V)\,.
\end{equation}
We note that the spatial finite element space $V_h^{p}$ is allowed to be
different on all subintervals $I_n$ which is natural in the context of a
discontinuous Galerkin approximation of the time variable and allows dynamic
mesh changes in time.
Due to the conformity of $V_h^{p}$ we get $\mathcal{V}_{\tau h}^{r,p}\subseteq
\mathcal{V}_{\tau}^{r}$.
Now, the fully discrete discontinuous in time scheme reads as follows:
\textit{Find $u_{\tau h} \in \mathcal{V}_{\tau h}^{r,p}$ such that}
\begin{equation}
\label{eq:11:A_tau_h_u_phi_eq_F_phi}
 A_{\tau}(u_{\tau h})(\varphi_{\tau h})
=
F_\tau(\varphi_{\tau h})
\quad \forall \varphi_{\tau h} \in \mathcal{V}_{\tau h}^{r,p}\,,
\end{equation}
\textit{with} $A_{\tau}(\cdot)(\cdot), a(\cdot)(\cdot)$ \textit{and} $F_\tau(\cdot)$
\textit{being defined in \eqref{eq:8:BilinearFormAtau},\eqref{eq:5:aBilinearform} and
\eqref{eq:9:LinearFormFtau}, respectively.}
In this work, we focus on convection-dominated problems with small diffusion
parameter $\varepsilon >0$. Then, the finite element approximation needs to be
stabilized in order to reduce spurious and non-physical oscillations of the
discrete solution arising close to layers.
Here, we apply the streamline upwind Petrov-Galerkin (SUPG) method
\cite{BruchhaeuserHB79,BruchhaeuserBH82}.
Existing convergence analyses in the natural norm of the underlying scheme
including local and global error bounds can be found, for instance, in~\cite[Ch. III.4.3]{BruchhaeuserRST08}.
The stabilized variant of the fully discrete discontinuous in time scheme
then reads as follows:
\textit{
Find $u_{\tau h} \in \mathcal{V}_{\tau h}^{r,p}$ such that
}
\begin{equation}
\label{eq:12:StabDGFully}
A_{S}(u_{\tau h})(\varphi_{\tau h}) =
F_{\tau}(\varphi_{\tau h}) \quad \forall \varphi_{\tau h}
\in \mathcal{V}_{\tau h}^{r,p}\,,
\end{equation}
\textit{
where the linear form $F_{\tau}(\cdot)$ is defined by
Eq.~\eqref{eq:9:LinearFormFtau} and the stabilized bilinear form
$A_{S}(\cdot)(\cdot)$ is given by}
\begin{equation}
\label{eq:13:ASutauhDG}
A_{S}(u_{\tau h})(\varphi_{\tau h}) \coloneqq
A_{\tau}(u_{\tau h})(\varphi_{\tau h})
+ S(u_{\tau h})(\varphi_{\tau h})\,.
\end{equation}
\textit{
Here, the SUPG stabilization bilinear form
}
$S(\cdot)(\cdot)$
\textit{
is defined by
}
\begin{equation}
\label{eq:14:StabDGSutauh}
\begin{array}{r@{\;}c@{\;}l@{\;}}
S(u_{\tau h})(\varphi_{\tau h}) & \coloneqq &
\displaystyle
\sum_{n=1}^{N}\int_{I_n}\sum_{K \in \mathcal{T}_{h, n}}
\delta_K\big(
r(u_{\tau h}),
\boldsymbol{b} \cdot \nabla \varphi_{\tau h}\big)_K \,\mathrm{d} t
\\[3ex]
& &
\displaystyle
+ \sum_{n=2}^{N}\sum\limits_{K\in\mathcal{T}_{h, n}}
\delta_K
\big(\left[u_{\tau h}\right]_{n-1},
\boldsymbol{b} \cdot \nabla \varphi_{\tau h,n-1}^+\big)_{K}
\\[3ex]
& &
+ \displaystyle\sum\limits_{K\in\mathcal{T}_{h, 1}}
\delta_K \big(u_{\tau h,0}^+ - u_0,
\boldsymbol{b} \cdot \nabla \varphi_{\tau h,0}^{+} \big)_{K}\,,
\end{array}
\end{equation}
\textit{
where the residual term $r(\cdot)$
is defined by}
\begin{equation}
\label{eq:15:rutauh}
r(u_{\tau h}) \coloneqq \partial_{t} u_{\tau h}
- \nabla\cdot\left(\varepsilon\nabla u_{\tau h}\right)
+ \boldsymbol{b} \cdot \nabla u_{\tau h}
+ \alpha u_{\tau h}
- f \,.
\end{equation}
\begin{remark}
\label{rem:2:StabParameter}
The proper choice of the stabilization parameter $\delta_K$ is an important
issue in the application of the SUPG approach; cf., e.g.,
\cite{BruchhaeuserJN11,BruchhaeuserJS08,BruchhaeuserJKN18}
and the discussion therein. For time-dependent CDR
problems an optimal error estimate for $\delta_K=\mathrm{O}(h)$ is derived
in~\cite{BruchhaeuserJN11}. In this work, we choose $\delta_{K}=\delta_0\cdot h_K$, where $\delta_0=0.1$ and $h_K=\sqrt[d]{|K|}$ is the cell diameter.
\end{remark}
Here and in the following we will use the notation cG($p$)-dG($r$) for a space-time finite elelemt Galerkin discretization with continuous pieceweise polynomials of degree $p$ in space and discontinuous piecewise polynomials of degree $r$ in time, cf., e.g., \cite{BruchhaeuserEEHJ96,BruchhaeuserBR12}.
\section{A Posteriori Error Representation}
\label{sec:error_rep}
In this section, we present an a posteriori error representation based on the
DWR method for the the stabilized problem \eqref{eq:12:StabDGFully}.
This representation is given in terms of a user-chosen goal functional $J \in\mathcal{C}^3(\mathcal{V},\mathbb{R})$, in general given as
\begin{equation} \label{def: Single Functional}
J(u)=\int_0^T J_t(u(t))\mathrm{d}t
+ J_T(u(T))\,,
\end{equation}
where $J_t$ and $J_T$ are three times differentiable functionals and each of
them may be zero; cf., e.g.,\cite{BruchhaeuserR17}. Since we focus
on anisotropic error estimation, we restrict this section to the main result and refer to~\cite[Ch.~4]{BruchhaeuserB22} for a detailed derivation in the context of CDR equations.

For the error representations in Theorem.~\ref{Thm:3.3}, we introduce the following Lagrangian
functionals
$\mathcal{L}: \mathcal{V}\times \mathcal{V} \rightarrow \mathbb{R}$,
$\mathcal{L}_\tau: \mathcal{V}_{\tau}^{r} \times \mathcal{V}_{\tau}^{r}
\rightarrow \mathbb{R}$, and
$\mathcal{L}_{\tau h}:
\mathcal{V}_{\tau h}^{r,p} \times \mathcal{V}_{\tau h}^{r,p}
\rightarrow \mathbb{R}$ by
\begin{subequations}
  \label{eq:3:3:Def_L_u_z_Def_L_tau_u_z_Def_L_tau_h_u_z}
  \begin{align}
    \label{eq:3:3:Def_L_u_z}
    \mathcal{L}(u,z) & \coloneqq  J(u)
                       + F(z)
                       - A(u)(z)\,,
    \\
    \label{eq:3:3:Def_L_tau_u_z}
    \mathcal{L}_{\tau}(u_\tau,z_\tau) & \coloneqq
                                        J(u_{\tau}) + F_\tau(z_{\tau})
                                        - A_{\tau}(u_{\tau})(z_{\tau})\,,
    \\
    \label{eq:3:3:Def_L_tau_h_u_z}
    \mathcal{L}_{\tau h}(u_{\tau h},z_{\tau h}) & \coloneqq
                                                  J(u_{\tau h})
                                                  + F_\tau (z_{\tau h})
                                                  - A_S(u_{\tau h})(z_{\tau h})=\mathcal{L}_{\tau}(u_{\tau h}),z_{\tau h}))-S(u_{\tau h}))(z_{\tau h}))\,.
  \end{align}
\end{subequations}
In~\eqref{eq:3:3:Def_L_u_z_Def_L_tau_u_z_Def_L_tau_h_u_z}, the Lagrange
multipliers $z$, $z_\tau$, and $z_{\tau h}$ are called adjoint variables in
contrast to the primal variables $u$, $u_\tau$, and $u_{\tau h}$;
cf.~\cite{BruchhaeuserBR12,BruchhaeuserBR01}.
The directional derivatives of the Lagrangian functionals
(G\^{a}teaux derivatives) with respect to their second argument yields the
primal
problems~\eqref{eq:2:WeakCDRsteady},~\eqref{eq:7:dGDiscTime},~\eqref{eq:12:StabDGFully},
while the directional derivative with respect to their first argument yields the
adjoint problems given by
\begin{subequations}
\label{eq:3:9:A_tau_prime_u_phi_z_eq_J_prime_u_phi_A_S_prime_u_phi_z_eq_J_prime_u_phi}
\begin{align}
  \label{eq:3:6:A_prime_u_phi_z_eq_J_prime_u_phi}
  A^{\prime}(u)(\varphi,z)
  &=
  J^{\prime}(u)(\varphi)
  \quad \forall \varphi \in \mathcal{V}\,,
  \\
  \label{eq:3:9:A_tau_prime_u_phi_z_eq_J_prime_u_phi}
  A_{\tau}^{\prime}(u_{\tau})(\varphi_{\tau},z_{\tau})
  & =
    J^{\prime}(u_{\tau})(\varphi_{\tau})
    \quad \forall \varphi_{\tau}\in \mathcal{V}_{\tau}^{r}\,,
  \\
  \label{eq:3:9:A_S_prime_u_phi_z_eq_J_prime_u_phi}
  A_{S}^{\prime}(u_{\tau h})(\varphi_{\tau h},z_{\tau h})
  & =
    J^{\prime}(u_{\tau h})(\varphi_{\tau h})
    \quad \forall \varphi_{\tau h}\in \mathcal{V}_{\tau h}^{r,p}\,,
\end{align}
\end{subequations}
where the definitions of $A^{\prime}(\cdot)(\cdot,\cdot)$, $A_{\tau}^{\prime}(\cdot)(\cdot, \cdot)$,
$A_{S}^{\prime}(\cdot)(\cdot, \cdot)$ are given in the appendix.
\begin{theorem}
\label{Thm:3.3}
Let $\{u,\,z\}\in \mathcal{V} \times \mathcal{V}$,
$\{u_{\tau},z_{\tau}\}
\in
\mathcal{V}_{\tau}^{r} \times \mathcal{V}_{\tau}^{r}$,
and
$\{u_{\tau h},z_{\tau h}\}
\in \mathcal{V}_{\tau h}^{r,p} \times \mathcal{V}_{\tau h}^{r,p}$
denote the stationary points of
$\mathcal{L}, \mathcal{L}_{\tau}$, and $\mathcal{L}_{\tau h}$
on the different levels of discretization, i.e.,
\begin{displaymath}
\begin{aligned}
\mathcal{L}^{\prime}(u,z)(\delta u, \delta z)
= \mathcal{L}_{\tau}^{\prime}(u,z)(\delta u, \delta z)
& = 0 \quad
\forall \{\delta u,\delta z\}\in \mathcal{V} \times \mathcal{V}\,,
\\
\mathcal{L}_{\tau}^{\prime}(u_{\tau},z_{\tau})
(\delta u_{\tau}, \delta z_{\tau})
& = 0
\quad \forall \{\delta u_{\tau},\delta z_{\tau}\}
\in \mathcal{V}_{\tau}^{r} \times \mathcal{V}_{\tau}^{r}\,,
\\
\mathcal{L}_{\tau h}^{\prime}(u_{\tau h},z_{\tau h})
(\delta u_{\tau h}, \delta z_{\tau h})
& = 0
\quad \forall \{\delta u_{\tau h},\delta z_{\tau h}\}
\in \mathcal{V}_{\tau h}^{r,p} \times \mathcal{V}_{\tau h}^{r,p}\,.
\end{aligned}
\end{displaymath}
Then, for the discretization errors in space and time we get
the representation formulas
\begin{subequations}
\label{eq:3:16}
\begin{align}
\label{eq:3:16a:J_u_minus_J_u_tau}
J(u)-J(u_{\tau}) & =
\frac{1}{2}\rho_{\tau}(u_{\tau})(z-\tilde{z}_{\tau})
+ \frac{1}{2}\rho_{\tau}^{\ast}(u_{\tau},z_{\tau})
(u-\tilde{u}_{\tau})
+ \mathcal{R}_{\tau}\,,
\\
\label{eq:3:16b:J_u_tau_minus_J_u_tau_h}
J(u_{\tau})-J(u_{\tau h}) & =
\frac{1}{2}\rho_{\tau}(u_{\tau h})(z_{\tau}-\tilde{z}_{\tau h})
+ \frac{1}{2}
\rho_{\tau}^{\ast}(u_{\tau h},z_{\tau h})
(u_{\tau}-\tilde{u}_{\tau h})
\\
\nonumber
& \qquad
+ \frac{1}{2} S(u_{\tau h})(\tilde{z}_{\tau h}+z_{\tau h})
+ \frac{1}{2} S^{\prime}(u_{\tau h})(\tilde{u}_{\tau h}-u_{\tau h},z_{\tau h})
+ \mathcal{R}_{h}\,,
\end{align}
\end{subequations}
where $\rho_{\tau}$ and $\rho_{\tau}^{\ast}$ are the primal and adjoint residuals based on the semi-discrete
in time schemes, respectively, given by
\begin{subequations}
  \begin{align}
    \label{eq:3:15:primal_dual_residuals}
    \rho_{\tau}(u)(\varphi)  &\coloneqq
    \mathcal{L}_{\tau,z}^{\prime}(u,z)(\varphi)=F_{\tau}(\varphi) - A_{\tau}(u)(\varphi)\,,\\
    \rho_{\tau}^{\ast}(u,z)(\varphi)
    &\coloneqq
    \mathcal{L}_{\tau,u}^{\prime}(u,z)(\varphi)= J^{\prime}(u)(\varphi) - A^{\prime}_{\tau}(u)(\varphi,z)
    \,.
  \end{align}
\end{subequations}
Here,
$\{\tilde{u}_{\tau},\tilde{z}_{\tau}\}\in \mathcal{V}_{\tau}^{r}
\times \mathcal{V}_{\tau}^{r}$,
and
$\{\tilde{u}_{\tau h},\tilde{z}_{\tau h}\} \in
\mathcal{V}_{\tau h}^{r,p} \times \mathcal{V}_{\tau h}^{r,p}$
can be chosen arbitrarily and $\mathcal{R}_{\tau}, \mathcal{R}_{h}$
are higher-order remainder terms with respect to the errors
$u-u_{\tau}, z-z_{\tau}$ and $u_{\tau}-u_{\tau h}, z_{\tau}-z_{\tau h}$, respectively.
\end{theorem}


\section{Anisotropic Error Estimators}
\label{sec:anisotropic-method}
The error representation formulas derived in Theorem~\ref{Thm:3.3} lead to a posteriori error estimators in space and time, which then serve as indicators for the adaptive mesh refinement process.
In order to obtain directional error indicators with regard to anisotropic mesh refinement, we split the spatial estimator into directional contributions using anisotropic interpolation and restriction operators.
To keep this section rather clear, we restrict ourselves to the presentation of the temporal and splitted spatial error estimators and refer to \cite{BruchhaeuserBBEMTW25} for a specific derivation and further details.
Throughout the remainder of this section, we restrict ourselves
to the two-dimensional case, i.\,e. $\Omega \subset \mathbb{R}^2$.

The temporal and directional in space ($i \in \{1,2\}$) error estimators are given by
\begin{subequations}\label{eq: implemented error estimator}
\begin{align}\label{eq: implemented error estimator-tau}
\eta_{\tau}&\coloneqq\frac{1}{2}\rho_{\tau}(u_{\tau h})(\mathrm{E}_{\tau}^{(r+1)}z_{\tau h}-z_{\tau h})
+ \frac{1}{2}\rho_{\tau}^{\ast}(u_{\tau h},z_{\tau h})(\mathrm{E}_{\tau}^{(r+1)}u_{\tau h}-u_{\tau h})\,,\\
\label{eq: implemented error estimator-h}
\eta_{h,i}&\coloneqq\frac{1}{2}\rho_{\tau}(u_{\tau h})(z_{\tau h}- R^{p}_iz_{\tau h})
+ \frac{1}{2}
\rho_{\tau}^{\ast}(u_{\tau h},\mathrm{R}_{h}^{p}z_{\tau h})
(\mathrm{I}_{2h,i}^{\left(2p\right)}u_{\tau h}-u_{\tau h})
\\
& \nonumber \qquad
+ \frac{1}{2} S(u_{\tau h})( R^{p}_iz_{\tau h})
+ \frac{1}{2} S^{\prime}(u_{\tau h})(\mathrm{I}_{2h,i}^{\left(2p\right)}u_{\tau h}-u_{\tau h},\mathrm{R}_{h}^{p}z_{\tau h})\,.
\end{align}
\end{subequations}
\begin{remark}\label{rem:negelecting}
Due to the splitting into directional contributions of the spatial error estimator, one obtains an additional anisotropic remainder term $\eta_{h,\mathbb{E}}$ that is of higher order and thus omitted in our numerical
examples. For a detailed derivation of this remainder term we refer to \cite[Sec.~4]{BruchhaeuserBBEMTW25}.
\end{remark}
The error estimators above have to be represented in a localized form in order to mark elements within the adaptive mesh refinement algorithm.
Here, we follow the classical way of the DWR philosophy, receiving the error representation on every mesh element by a cell-wise integration by parts, cf.~\cite{BruchhaeuserBR01,BruchhaeuserBR03}.
The decomposition into elementwise contributions of the anisotropic error indicators in space and time are given by
\begin{equation}
\label{eq:error_indicators}
\eta_{\tau}=\sum_{n=1}^{N}\sum_{K\in\mathcal{T}_h^n}\eta_{\tau}^{K,n}\,,
\qquad
\eta_{h,i}=\sum_{i=1}^d\sum_{n=1}^{N} \sum_{K\in\mathcal{T}_h^n} \eta_{h,i}^{K,n}\,.
\end{equation}

\subsection{Practical Aspects}
\label{subsection:practical_aspects}
\begin{figure}[htb]
  \begin{minipage}{0.245\textwidth}
    \scalebox{0.6}{\begin{tikzpicture}
\begin{axis}[
xlabel={$x$},                  
ylabel={$y$},               
hide z axis,
grid=both,                    
samples=3,                   
domain=0:2,                   
y domain=0:2,                 
xtick={0,1,2},
xticklabels={},
yticklabels={},
xlabel={$x_1$},
ylabel={$x_2$},
colormap/jet,
mesh/ordering=x varies,
mesh/cols=13
]
\addplot3[
surf,                         
shader=interp,                  
opacity=0.3                   
]
{\functionforplot{x}{y}};  

\addplot3[
mesh,
draw=black,
samples=3, 
]
{\functionforplot{x}{y}};  
\end{axis}
\end{tikzpicture} }
    \subcaption{$v_h$.}
  \end{minipage} \hfill%
  \begin{minipage}{0.245\textwidth}
    \scalebox{0.6}{\begin{tikzpicture}
\begin{axis}[
xlabel={$x$},                 
ylabel={$y$},     
hide z axis,           
grid=both,                    
samples=30,                   
samples y=30,
domain=0:2,                   
y domain=0:2,                 
xtick={0,1,2},
xticklabels={},
yticklabels={},
xlabel={$x_1$},
ylabel={$x_2$},
colormap/jet,
mesh/ordering=x varies,
mesh/cols=13
]

\addplot3[
surf,                         
shader=interp,  
samples=3, 	              
samples y=30,
opacity=0.3                   
]
{\functionforplot{x}{y}};  
\addplot3 [domain=0:2,samples y=1,color=red,ultra thick] (1,x,\functionforplot{1}{x});
\addplot3 [domain=0:2,samples y=1,color=red,ultra thick] (0,x,\functionforplot{0}{x});
\addplot3 [domain=0:2,samples y=1,color=red,ultra thick] (2,x,\functionforplot{2}{x});
\addplot3 [domain=0:2,samples =3,samples y=1] (x,1,\functionforplot{x}{1});
\addplot3 [domain=0:2,samples =3,samples y=1] (x,0,\functionforplot{x}{0});
\addplot3 [domain=0:2,samples =3,samples y=1] (x,2,\functionforplot{x}{2});
\end{axis}

\end{tikzpicture} }
    \subcaption{$\mathrm{I}_{2h,1}^{\left(2p\right)}v_h$.}
  \end{minipage} \hfill%
  \begin{minipage}{0.245\textwidth}
    \scalebox{0.6}{\begin{tikzpicture}
\begin{axis}[
xlabel={$x$},                 
ylabel={$y$},     
hide z axis,           
grid=both,                    
samples=30,                   
samples y=3,
domain=0:2,                   
y domain=0:2,                 
xtick={0,1,2},
xticklabels={},
yticklabels={},
xlabel={$x_1$},
ylabel={$x_2$},
colormap/jet,
mesh/ordering=x varies,
mesh/cols=13
]

\addplot3[
surf,                         
shader=interp,  	              
samples y=3,
opacity=0.3                   
]
{\functionforplot{x}{y}};  
\addplot3 [domain=0:2,samples =3,samples y=1] (1,x,\functionforplot{1}{x});
\addplot3 [domain=0:2,samples =3,samples y=1] (0,x,\functionforplot{0}{x});
\addplot3 [domain=0:2,samples =3,samples y=1] (2,x,\functionforplot{2}{x});
\addplot3 [domain=0:2,samples y=1,color=red,ultra thick] (x,1,\functionforplot{x}{1});
\addplot3 [domain=0:2,samples y=1,color=red,ultra thick] (x,0,\functionforplot{x}{0});
\addplot3 [domain=0:2,samples y=1,color=red,ultra thick] (x,2,\functionforplot{x}{2});
\end{axis}

\end{tikzpicture}}
    \subcaption{$\mathrm{I}_{2h,2}^{\left(2p\right)}v_h$.}
  \end{minipage}\hfill%
  \begin{minipage}{0.245\textwidth}
    \scalebox{0.6}{\begin{tikzpicture}
\begin{axis}[
xlabel={$x$},                  
ylabel={$y$},               
hide z axis,
grid=both,                    
samples=30,                   
domain=0:2,                   
y domain=0:2,                 
xtick={0,1,2},
xticklabels={},
yticklabels={},
xlabel={$x_1$},
ylabel={$x_2$},
colormap/jet,
mesh/ordering=x varies,
mesh/cols=13
]
\addplot3[
surf,                         
shader=interp,                  
opacity=0.3                   
]
{\functionforplot{x}{y}};  

\addplot3 [domain=0:2,samples y=1] (1,x,\functionforplot{1}{x});
\addplot3 [domain=0:2,samples y=1] (0,x,\functionforplot{0}{x});
\addplot3 [domain=0:2,samples y=1] (2,x,\functionforplot{2}{x});
\addplot3 [domain=0:2,samples y=1] (x,1,\functionforplot{x}{1});
\addplot3 [domain=0:2,samples y=1] (x,0,\functionforplot{x}{0});
\addplot3 [domain=0:2,samples y=1] (x,2,\functionforplot{x}{2});
\end{axis}
\end{tikzpicture}}
    \subcaption{$\mathrm{I}_{2h}^{\left(2p\right)}v_h$.}
  \end{minipage}%
  \caption{The action of different interpolation operators on a function $v_h$ on the path $K_{2h}$ and $p=1$.\label{Fig: Action AnIso Int}}
\end{figure}
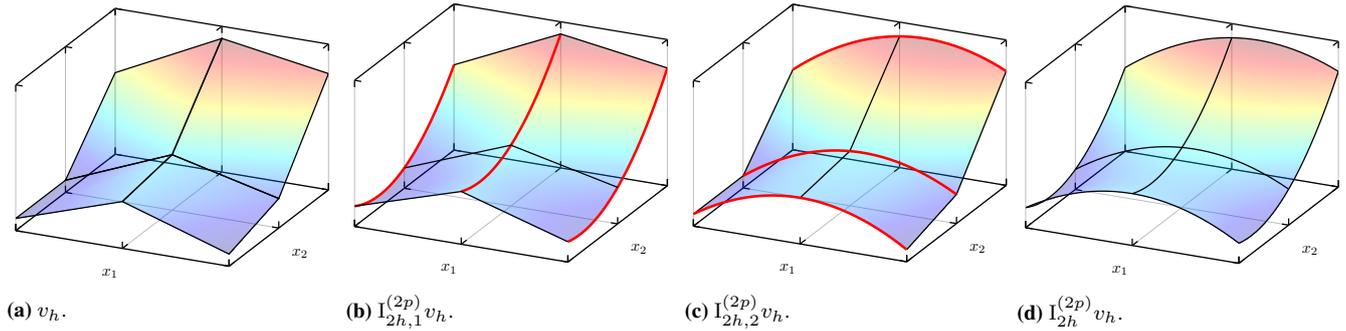

In this subsection, we explain the single steps needed to obtain the computable versions of the error estimators given by Eq.~\eqref{eq: implemented error estimator} from the error representaion formulas derived in Theorem~\ref{Thm:3.3}.
First, we neglect the remainder terms  $\mathcal{R}_{\tau}$ and $\mathcal{R}_{h}$ in~\eqref{eq:3:16} as they are of higher order.
Second, the unknown solutions in~\eqref{eq:3:16} are replaced by means of the approximated fully discrete solutions $u_{\tau h}\in\mathcal{V}_{\tau h}^{r,p}$ and $z_{\tau h}\in \mathcal{V}_{\tau h}^{r,q}$, with $q=2p$, respectively.
Third, the temporal and spatial weights have to be approximated in a suitable
way, where we refer to \cite{BruchhaeuserBBEMTW25} for the detailed definitions of the following (anisotropic) interpolation and restriction operators as well as to \cite{BruchhaeuserBR01,BruchhaeuserBR03} for a short review of possible techniques in general.
The following choice of approximations for the temporal and spatial weights are based on our previous numerical comparisons in~\cite{BruchhaeuserBB24,BruchhaeuserBSB19}:
\paragraph{Temporal Weights.} Approximate the temporal weights $u-\tilde{u}_{\tau}$ and $z-\tilde{z}_{\tau}$,
respectively, by means of a higher-order reconstruction using Gauss
quadrature points given by
\begin{displaymath}
u-\tilde{u}_{\tau} \approx \mathrm{E}_{\tau}^{(r+1)}u_{\tau h}-u_{\tau h}\,,
\qquad
z-\tilde{z}_{\tau} \approx \mathrm{E}_{\tau}^{(r+1)}z_{\tau h}-z_{\tau h}\,,
\end{displaymath}
using a reconstruction in time operator $\mathrm{E}_{\tau}^{(r+1)}$
thats lifts the solution to a piecewise
polynomial of degree ($r+1$) in time, cf.~\cite{BruchhaeuserBB24,BruchhaeuserB22} for further details.
\paragraph{Spatial Weights.} Approximate the spatial weights $u_{\tau}-\tilde{u}_{\tau h}$ and $z_{\tau}-\tilde{z}_{\tau h}$,
respectively, by means of anisotropic interpolation and restriction operations given by
\begin{equation}
\label{eq:spatial_weights}
u_{\tau}-\tilde{u}_{\tau h}  \approx  \mathrm{I}_{2h,i}^{(2p)}u_{\tau h}-u_{\tau h}\,,
\qquad
z_{\tau}-\tilde{z}_{\tau h}  \approx  z_{\tau h}-\mathrm{R}_{i}^{p}z_{\tau h}\,,
\qquad
\tilde{z}_{\tau h} \approx  \mathrm{R}_{h}^{p}z_{\tau h}\,.
\end{equation}
Here, $R^p_h:\mathcal{V}_{h}^{2p} \to \mathcal{V}_{h}^{p}$ denotes a restriction in space operator that acts on a spatial cell and restricts the solution to a polynomial of degree $p<q$ on the corresponding reference cell (cf.~\cite[Def. 3.6]{BruchhaeuserBBEMTW25}).
Moreover, $R^p_i:{V}_{ h}^{2p} \to {V}_{ h}^{2p}$ denotes the restriction in space operator in $i$-th direction that acts on a spatial cell and restricts the solution to a polynomial of degree $p<q$ in $i$-th direction on the corresponding reference cell (cf.~\cite[Def. 4.3]{BruchhaeuserBBEMTW25}).
Finally, $\mathrm{I}_{2h,i}^{\left(2p\right)}\colon V_h^{p}(K_{2h}^P) \to V_h^{2p}(K_{2h})$, defined as
\begin{equation}
\label{eq:I2hi2p}
\mathrm{I}_{2h,i}^{\left(2p\right)}\coloneqq R^{p}_i\circ\mathrm{I}_{2h}^{\left(2p\right)}
\end{equation}
denotes a patch-wise higher order interpolation in space operator in $i$-th direction that acts on a patched cell of size $2h$ and lifts the solution to a piecewise polynomial of degree $2p$ in $i$-th direction on the reference cell corresponding to the patched cell of width $2h$ (cf.~Fig.~\ref{Fig: Action AnIso Int} and \cite[Def. 4.9]{BruchhaeuserBBEMTW25}). In \eqref{eq:I2hi2p},
$\mathrm{I}_{2h}^{\left(2p\right)}$ denotes the patch-wise higher order interpolation in space operator, cf., e.g., \cite[Def.~3.4]{BruchhaeuserBBEMTW25} or \cite[Ch.~4]{BruchhaeuserBR03}. Further, the operator $R_i^p$, which was originally defined for elements $K \in \mathcal{T}_h$, is here extended to act on elements $K_{2h} \in \mathcal{T}_{2h}$.
The action of directional interpolation operators, exemparily for the case $p=1$, is visualized in Figure~\ref{Fig: Action AnIso Int}.
For more details we refer to \cite[Sec. 4]{BruchhaeuserBBEMTW25}.


\subsection{Multi-goal oriented error control}
\label{sec:4:2:multi-goal}
In many applications, multiple quantities of interest (goal functionals) are considered.  Suppose we have $N \in \mathbb{N}$, where we denote $\boldsymbol{J} \coloneqq \left( J_1, J_2, J_3, \ldots, J_N \right)$, with each $J_i$ as defined in \eqref{def: Single Functional}. From a general perspective, one might question if uniform mesh refinement might be more efficient for several goal functionals. However, the works \cite{BruchhaeuserHH03,Pardo10,EnLaWi20} have demonstrated that adaptive refinement still is more effective.  While it is possible to individually approximate the error estimator for each goal functional, this approach requires solving the adjoint problem $N$ times. Inspired by the works of \cite{BruchhaeuserHH03, BruchhaeuserH08, Pardo10, EnLaWi20}, we combine these goal functionals into a single composite functional $J_w$ defined as

\begin{equation}
  J_w(v) \coloneqq \sum_{i=1}^{N} w_i J_i(v) \qquad \forall v \in V,
\end{equation}
where $w_i \in \mathbb{R}$ are weights.
As discussed in \cite{EndtmayerWick17, EndLanWic19, EnLaWi20}, the choice of the sign for these weights is critical to avoid error cancellation effects. It turns out that $w_i(J_i(u) - J_i(u_h))>0$ is sufficient to satisfy this. However, since $J_i(u)$ is not known, it can only be approximated. In the studies \cite{EndLanWic19, EnLaWi20}, $J(u)$ is approximated by $J(u^{(2)})$, where $u^{(2)}$ represents a higher order approximation of the solution.
Once the weights $w_i$ are chosen, the resulting adjoint problems for $J_w$ are given by
\begin{subequations}
  \begin{align*}
    A^{\prime}(u)(\varphi,z_w)
    &=
    J^{\prime}_w(u)(\varphi)
    \quad \forall \varphi \in \mathcal{V}\,,
    \\
    A_{\tau}^{\prime}(u_{\tau})(\varphi_{\tau},z_{\tau w})
    & =
    J^{\prime}_w(u_{\tau})(\varphi_{\tau})
    \quad \forall \varphi_{\tau}\in \mathcal{V}_{\tau}^{r}\,,
    \\
    A_{S}^{\prime}(u_{\tau h})(\varphi_{\tau h},z_{\tau h w})
    & =
    J^{\prime}_w(u_{\tau h})(\varphi_{\tau h})
    \quad \forall \varphi_{\tau h}\in \mathcal{V}_{\tau h}^{r,p}\,.
  \end{align*}
\end{subequations}
Finally, we replace the $(J, z,z_{\tau},z_{\tau h})$ by $(J_w, z_w,,z_{\tau w},z_{\tau h w})$ in \eqref{eq: implemented error estimator-tau} and \eqref{eq: implemented error estimator-h}. For simplicity and to enhance readability, we retain the notation \((J, z, z_{\tau}, z_{\tau h})\) in further correspondence.

%
\subsection{Goal-Oriented Anisotropic Adaptive Algorithm}
\label{subsection:algorithm}
In this subsection, we present a survey-like version of our anisotropic adaptive solution algorithm using tensor-product space-time slabs $\mathcal{Q}_n$, where we refer to \cite{BruchhaeuserB22,BruchhaeuserBBEMTW25,BruchhaeuserKBB19,endtmayer_goal-oriented_2024,BruchhaeuserELRSW24,RoThiKoeWi23}
for algorithmic details and implementation aspects, especially regarding space-time finite elements.

To compute the error indicators required in Algorithm~\ref{alg:AnisotropicAMR}, we first solve the primal problem forward
in time, followed by the adjoint problem backward in time. Using both primal  and adjoint solutions, we calculate the local error indicators~\eqref{eq:error_indicators}, perform
marking and refinement according to Algorithm~\ref{alg:AnisotropicAMR}, and  repeat the space-time solution procedure iteratively until the desired
accuracy is reached.
\begin{algorithm}[htb]
\caption{\label{alg:AnisotropicAMR}Anisotropic Mesh Adaptation: mark \& refine}
\begin{algorithmic}[1]
\STATE{Set $\ell=1$ and generate initial space-time slabs
$\mathcal{Q}_n^1=\mathcal{T}_{h}^1\times I_n^1\,, n=1,\dots,N^1$.}
\STATE{Compute $u_{\tau h}^\ell$ by solving the \textbf{stabilized primal}
problem \eqref{eq:12:StabDGFully}.}\label{alg:primal}
\STATE{Compute $z_{\tau h}^\ell$ by solving the \textbf{stabilized adjoint}
problem \eqref{eq:3:9:A_S_prime_u_phi_z_eq_J_prime_u_phi}.}
\REQUIRE{$\displaystyle\eta^{K,n,\ell}_{h,i}$,
$\eta^{K,n,\ell}_{\tau}$ (cf.~\eqref{eq:error_indicators}) for all $K\in
\mathcal T_h^\ell$, for all $I_n^\ell\in \mathcal T_{\tau}^\ell$, $\theta_h^\ell$, $\theta_\tau^\ell$}
\STATE{\textbf{Calculate} for all $\displaystyle K\in \mathcal T_{h}^\ell,\:i=1,\dots,d$:
$\displaystyle\eta^{K,\ell}_{h,i}=\displaystyle\sum_{I_n^\ell\in \mathcal T_{\tau}^\ell}\eta^{K,n,\ell}_{h,i}$.}\label{alg:calc-eta-khi}
\STATE{\textbf{Calculate} for all $\displaystyle I_n^\ell\in \mathcal T_{\tau}^\ell$:        $\displaystyle\eta^{n,\ell}_{\tau}=\sum_{K\in \mathcal T_{h}^\ell}\eta^{K,n,\ell}_{\tau}$.}\label{alg:calc-eta-nt}
\FOR{$i=1,\dots,d$}
\STATE{\textbf{Mark} $\displaystyle\theta_h^\ell \lvert\mathcal T_{h}^\ell\rvert$
cells $\displaystyle K\in \mathcal
T_h^\ell$ with the largest refinement criteria $\displaystyle \eta^{K,\ell}_{h,i}$ for \emph{refinement in direction} $i$.}\label{alg:mark-aniso}
\ENDFOR{}
\STATE{\textbf{Mark} $\displaystyle\theta_\tau^\ell \lvert\mathcal T_{\tau}^\ell\rvert$ subintervals $\displaystyle I_n^\ell\in \mathcal
T_\tau^\ell$ with the largest refinement criteria $\displaystyle \eta^{n,\ell}_{\tau}$ for refinement.}\label{alg:mark-time}
\STATE{\textbf{Refine} in space according to marking in line~\ref{alg:mark-aniso}}
\STATE{\textbf{Refine} in time according to marking in line~\ref{alg:mark-time}}
\STATE{Increase $\ell$ to $\ell +1$ and return to line~\ref{alg:primal}}
\end{algorithmic}
\end{algorithm}

\begin{remark}[Anisotropic Mesh Adaptation]
~\\[-2\topsep]
\begin{itemize}\itemsep1pt \parskip0pt \parsep0pt
\item The spatial mesh is fixed for all time intervals, the local error
indicators \(\eta^{K,n}_{h,i}\), computed for each cell \(K\) and direction
\(i\), are summed over time intervals to yield the total directional error
\(\eta^K_{h,i}\) (Algorithm~\ref{alg:AnisotropicAMR},
line~\ref{alg:calc-eta-khi}).
The temporal error \(\eta^{n}_{\tau}\) results from summing local    contributions across all spatial elements (Algorithm~\ref{alg:AnisotropicAMR}, line~\ref{alg:calc-eta-nt}).
\item Marking for refinement is performed separately in each spatial direction
(Algorithm~\ref{alg:AnisotropicAMR}, line~\ref{alg:mark-aniso}). A spatial
element \(K\) is marked for refinement in direction \(i\) if its directional  indicator \(\eta^K_{h,i}\) is among the top \(\theta_h |\mathcal T_{h}|\)
values. Evaluating directional errors individually naturally yields anisotropic
refinement when one direction dominates, and isotropic refinement when
errors are balanced.
\item Adaptive refinement on quadrilateral or hexahedral meshes leads to
interdependent hanging nodes. To resolve these mutual dependencies, we refine the coarser element at the end of the hanging-node chain in the appropriate
direction (i.\,e.\ anisotropically).
\end{itemize}
\end{remark}

\section{Numerical Examples}
\label{sec:5:numerical_examples}
In this section, we investigate our anisotropic space-time adaptive algorithm with regard to multiple goals.
For simplicity, we set the weights to $w_i = 1$.
Thereby, we study accuracy, efficiency as well as stability properties by means of the following numerical examples.
\begin{enumerate}
\item \emph{Moving Hump with Circular Layer:}
For this benchmark case an exact solution is available, which is characterized by a circular, interior layer changing its height in the course of time.
\item \emph{Nonstationary Hemker Problem with circular obstacle:}
The Hemker problem serves as a more challenging test case. It can be interpreted as a model of a convection-dominated heat transfer from a hot column. The solution to the Hemker problem is characterized by a boundary layer at the cylinder's surface and two interior layers located downstream of the cylinder.
\end{enumerate}
To measure the accuracy of the error estimator we will study the \emph{effectivity index}
\begin{equation}
I_{\text{eff}} = \left\vert\frac{\eta_{h,\,x}+\eta_{h,\,y}+\eta_{\tau}}{\boldsymbol J(u) - \boldsymbol J(u_{\tau h})}\right\vert\,,
\end{equation}
defined as the ratio of the sum of the estimated errors over the exact error. Here, we denote
$\eta_{h}=\eta_{h,\,x}+\eta_{h,\,y}$ and
$\eta_{\tau h}=\eta_{\tau}+\eta_{h}$, where $\eta_{h,\,x} \coloneqq \eta_{h,1}$ and  $\eta_{h,\,y} \coloneqq \eta_{h,2}$.
Moreover, as an indicator for the anisotropy we consider the maximum
\emph{aspect ratio}, given by
\begin{equation}
\operatorname{ar}_{\max}\coloneqq \max_{K \in \mathcal T_h} \max_{\boldsymbol x \in
Q_K}\frac{\lambda_{\max}}{\lambda_{\min}}\,,
\end{equation}
where $Q_K$ is the set of quadrature points on $N_t$ and $\lambda_{\min},
\lambda_{\max}$ are the minimal and maximal eigenvalues of \(\nabla
\boldsymbol T_K(\boldsymbol x)\), respectively. The SUPG stabilization paramater $\delta_{K}(\cdot)$  is defined as
\[
\delta_K = \delta_0\,h_K\,,\;0.1\leq \delta_0\leq1\,,
\]
where $h_K=\sqrt[d]{|K|}$ is the cell diameter of the mesh cell.
The total number of DoFs  is denoted by $N_{\text{tot}}$, whereas the spatial DoFs and temporal DoFs are denoted by $N_{\text{space}}$ and $N_{\text{time}}$, respectively.
The implementation is based on the \texttt{deal.II}
finite element library~\cite{africa_dealii_2024}. The tests are run on a
single node with 2 Intel Xeon Platinum 8360Y CPUs and \SI{1024}{\giga\byte}RAM
of the HPC cluster HSUper at Helmut Schmidt University.\@

\subsection{Moving Hump with Circular Layer}
\label{sec:MovingHump}
The first example is a benchmark for convection-dominated problems that goes back to \cite{BruchhaeuserJS08}.
The (analytical) solution is characterized by a hump changing its height in the course of the time, given by
\begin{equation}
\label{eq:10:uexactMH}
\begin{array}{r@{\;}c@{\;}l@{\;}}
u(\boldsymbol{x}, t) :=
\frac{16}{\pi}\sin(\pi t)x_1(1-x_1)x_2(1-x_2)
\cdot
\big( \frac{1}{2} +
\arctan\big[2\varepsilon^{-\frac{1}{2}}\big(
r_0^2-(x_1-x_1^0)^2-(x_2-x_2^0)^2\big)
\big]
\big)
\end{array}
\end{equation}
The problem is defined on $\Omega\times I := (0,1)^2\times (0,0.5]$ with the
scalars $r_0=\frac{1}{4}$ and $x_1^0=x_2^0=\frac{1}{2}$.
We choose $\boldsymbol{b} = (2,3)^{\top}$ and $\alpha = 1$.
The goal quantity is chosen to control the error within two spatial control
points
\begin{equation}
  \label{}
  \boldsymbol x_{e_1}=\left(\frac{1}{2}-\frac{1}{\sqrt{2}}r_0,\,\frac{1}{2}-\frac{1}{\sqrt{2}}r_0\right)\quad\text{and}\quad \boldsymbol x_{e_2}=\left(\frac{1}{2}+\frac{1}{\sqrt{2}}r_0,\,\frac{1}{2}+\frac{1}{\sqrt{2}}r_0\right)\,,
\end{equation}
respectively, which are located within the circular layer.
Corresponding to Sec.~\ref{sec:4:2:multi-goal}, we use the following multi-goal
\begin{equation}\label{eq:point-goal}
\boldsymbol J(u)=\big(J_1(u),\,J_2(u)\big)\,,\textnormal{ with } \quad
J_1(u)=u(\boldsymbol x_{e_1},\,T)\,, \quad  J_2(u)=u(\boldsymbol x_{e_2},\,T)\,.
\end{equation}
Here, the weights are set to one ($w_i=1\,,i\in\{1,2\}$) such that
\[J_w(u)=\displaystyle\sum_{i=1}^2J_i(u)\,.\] Furthermore, we regularize $J_i$, $i\in\{1,2\}$, by a Dirac delta function $\delta_{r,\boldsymbol c}(\boldsymbol x)=\alpha\operatorname{e}^{(1-1/(1-r^2/s^2))}$, where $r=\lVert \boldsymbol x-\boldsymbol c\rVert$, $s>0$ is the cutoff radius and $\alpha$ the scaling factor, such that $\delta_{r,\boldsymbol c}$ integrates to $1$.
This test case is solved using a cG($1$)-dG($0$) discretization for the primal problem and a cG($2$)-dG($0$) discretization for the adjoint problem.

In Table~\ref{tab:movinghump}, we present the development of the total discretization error $\boldsymbol J(u)-\boldsymbol J(u_{\tau h})$ for the multi-goal functional \eqref{eq:point-goal} as well as the error indicators,
effectivity indices and aspect ratios during an anisotropic adaptive refinement process.
With regard to accuracy, we observe a very good approximation quality of the actual discretization error, as visible by effectivity indices very close to one in the course of the adaptive run.
Moreover, a good equilibration of the directional error indicators in space $\eta_{h, x}$ and $\eta_{h, y}$ is obtained in the final loops (cf. columns seven and eight of Table \ref{tab:movinghump}).
This equilibration is essential for the efficiency of the anisotropic adaptive algorithm.

\begin{figure}[h]
  \captionsetup[subfigure]{justification=centering}
\centering
\subfloat[]{
  \centering
\begin{minipage}{.49\linewidth}
\centering
\includegraphics[width=6.65cm]{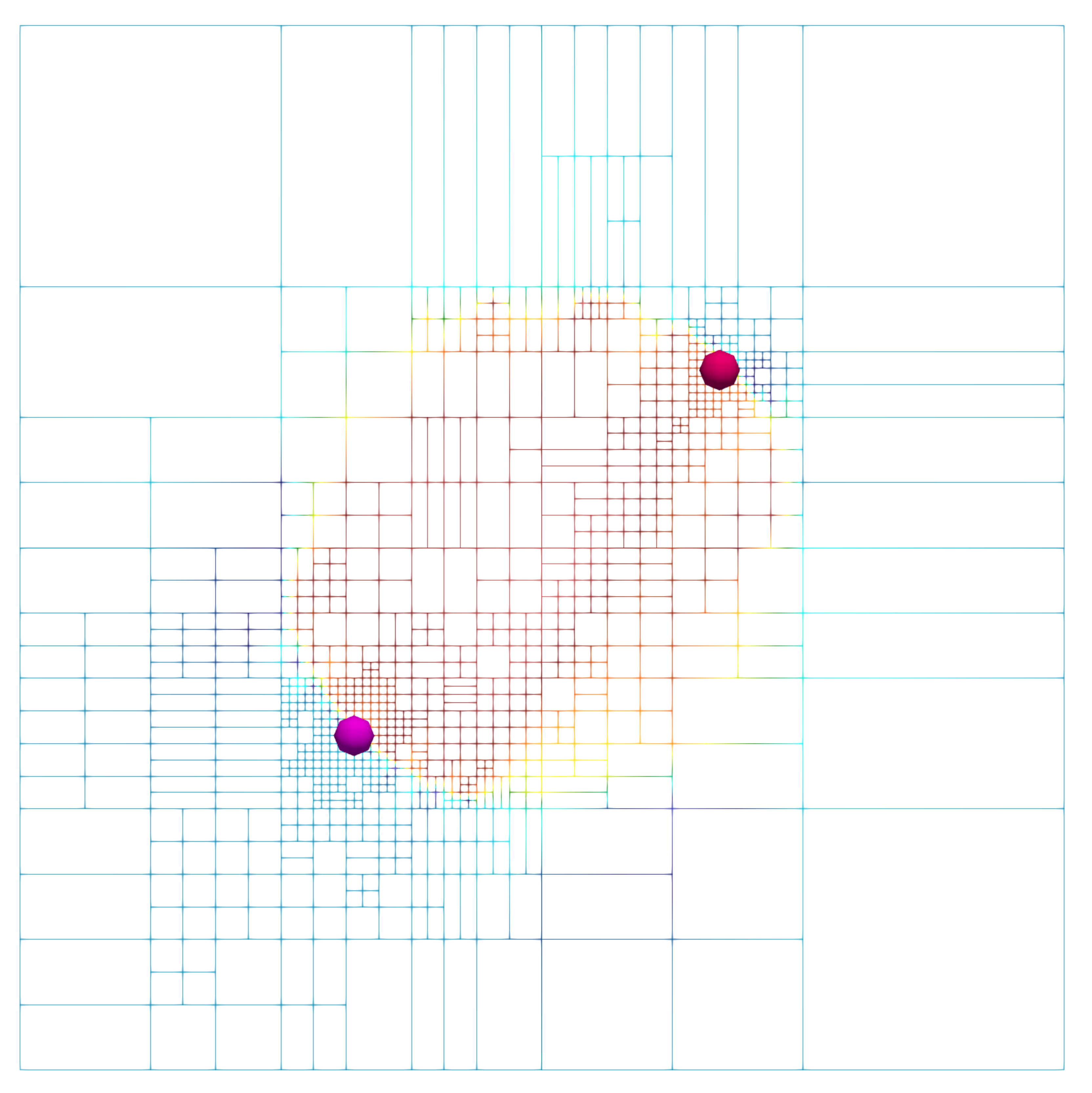}
\end{minipage}
\label{BruchhaeuserFig:0}
}
\subfloat[]{
  \centering
\begin{minipage}{.49\linewidth}
\centering
\includegraphics[width=6.65cm]{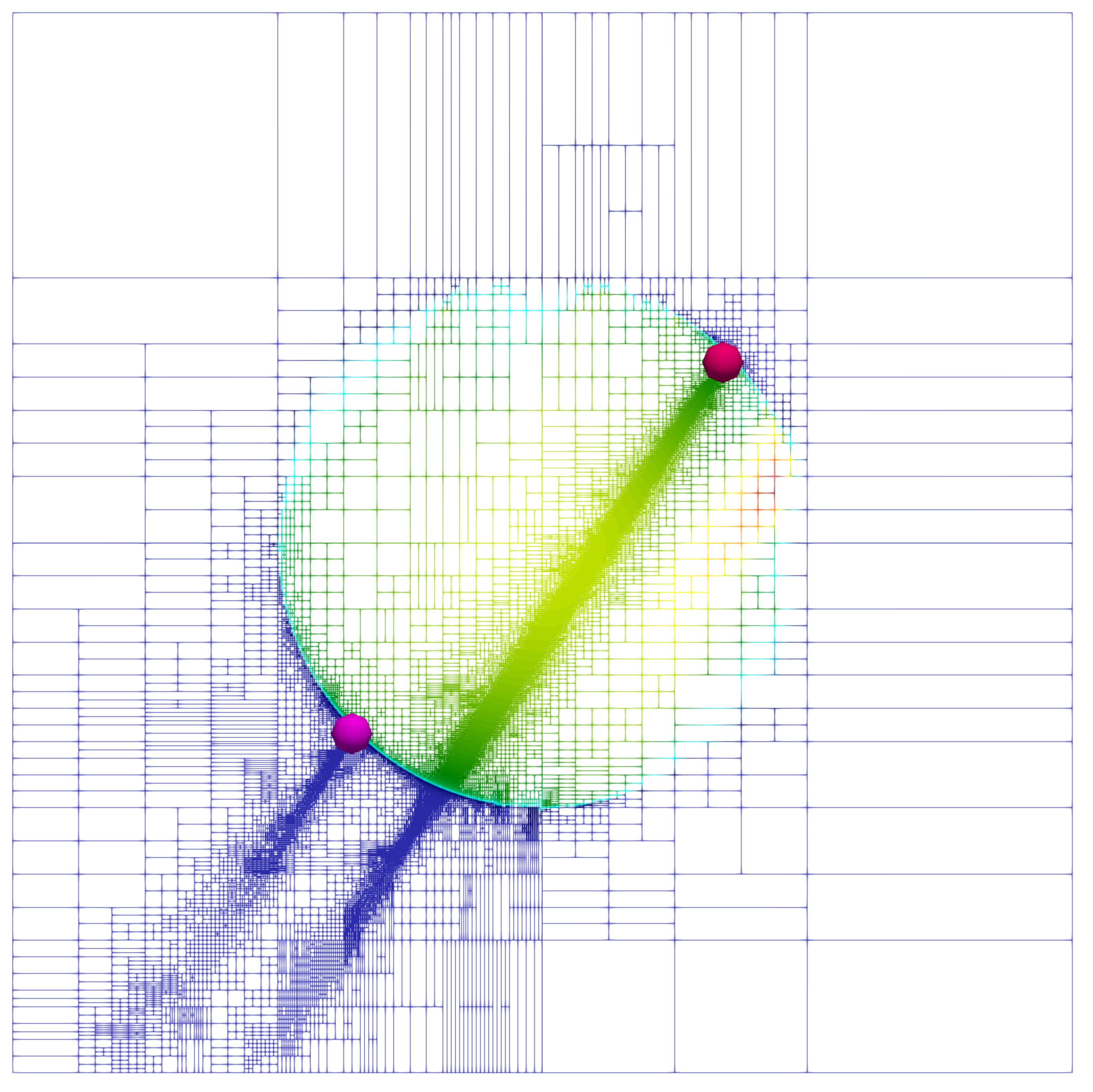}
\end{minipage}
\label{BruchhaeuserFig:25}
}
\caption{Anisotropic adaptive spatial meshes after
6 (a) and 13 (b) DWR loops at final time point
$T=0.5$ for $\varepsilon=10^{-6}$.}
\label{BruchhaeuserFig:SolutionProfiles}
\end{figure}

Finally, in Fig.~\ref{BruchhaeuserFig:SolutionProfiles} we visualize the anisotropic, adaptively refined spatial meshes for two different DWR loops, one of them at the outset and the other at the end.
We note that at the beginning the adaptive refinement is mainly located close to the two control points $\boldsymbol x_{e_1}$ and $\boldsymbol x_{e_2}$, whereas at the end the mesh is additionally refined along those cells that affect the multi-goal error by means of transport in the direction of the convection field $\boldsymbol b$. Moreover, mesh cells without strong impact on the solution close to the specific control points remain coarse. This behavor is desirably since the multi-goal aims to control the error only within the chosen points of control, leading to a high economical mesh.
\begin{table}[htb]\setlength{\tabcolsep}{3.5pt}
\caption{\label{tab:movinghump}
\emph{Anisotropic} adaptive refinement including effectivity indices, error indicators and maximum aspect ratios for goal quantity \eqref{eq:point-goal}, $\varepsilon=10^{-6}$, $\delta_0=1.0$ for Example~\ref{sec:MovingHump}.
}
\begin{center}\footnotesize
\begin{tabular}{crrrrrrrrrrr}
\toprule
$\ell$ & $N_{\text{space}}$ & $N_{t}$ & $N_{\text{tot}}$ &
{$\boldsymbol J(u)-\boldsymbol J(u_{\tau h})$}
& $\eta_{h,x}$  & $\eta_{h,y}$  & $\eta_h$ & $\eta_\tau$ &
$\eta_{\tau h}$ & $I_{\text{eff}}$ & $\operatorname{ar}_{\max}$
\\ \midrule
  \num{1 }& \num{ 25}     &\num{ 5 } & \num{125}     & ${8.416e\!-\!1}$&${ 1.703e\!-\!1 }$&${-3.253e\!-\!1 }$&${-4.957e\!-\!1 }$&${-1.263e\!-\!1 }$&${-6.220e\!-\!1 }$&${ 0.74 }$& \num{1  }\\
  \num{2 }& \num{49}      &\num{ 6 } & \num{294}     & ${6.738e\!-\!1}$&${ 2.055e\!-\!1 }$&${ 2.943e\!-\!1 }$&${ 4.998e\!-\!1 }$&${-5.264e\!-\!2 }$&${ 4.472e\!-\!1 }$&${ 0.66 }$& \num{1  }\\
  \num{3 }& \num{105}     &\num{ 7 } & \num{735}     & ${5.474e\!-\!1}$&${ 7.780e\!-\!1 }$&${ 1.236e\!+\!0 }$&${ 2.015e\!+\!0 }$&${ 2.479e\!-\!1 }$&${ 2.262e\!+\!0 }$&${ 4.13 }$& \num{2  }\\
  \num{4 }& \num{230}     &\num{ 8 } & \num{1840}    & ${1.857e\!-\!1}$&${ 2.030e\!+\!0 }$&${ 1.857e\!+\!0 }$&${ 3.888e\!+\!0 }$&${-4.390e\!-\!2 }$&${ 3.844e\!+\!0 }$&${20.70 }$& \num{4  }\\
  \num{5 }& \num{494}     &\num{ 9 } & \num{4446}    & ${3.618e\!-\!1}$&${ 3.402e\!+\!0 }$&${ 2.639e\!+\!0 }$&${ 6.042e\!+\!0 }$&${-4.865e\!-\!2 }$&${ 5.993e\!+\!0 }$&${16.57 }$& \num{8  }\\
  \num{6 }& \num{1064}    &\num{ 10} & \num{10640}   & ${1.783e\!-\!1}$&${-3.999e\!-\!1 }$&${-1.497e\!-\!2 }$&${-4.149e\!-\!1 }$&${ 4.986e\!-\!2 }$&${-3.651e\!-\!1 }$&${ 2.05 }$& \num{8  }\\
  \num{7 }& \num{2310}    &\num{ 11} & \num{25410}   & ${8.828e\!-\!2}$&${-1.039e\!-\!1 }$&${-1.279e\!-\!1 }$&${-2.319e\!-\!1 }$&${ 2.951e\!-\!2 }$&${-2.023e\!-\!1 }$&${ 2.29 }$& \num{16 }\\
  \num{8 }& \num{5001}    &\num{ 12} & \num{60012}   & ${4.423e\!-\!2}$&${ 1.600e\!-\!1 }$&${ 1.266e\!-\!1 }$&${ 2.867e\!-\!1 }$&${ 2.225e\!-\!2 }$&${ 3.089e\!-\!1 }$&${ 6.99 }$& \num{16 }\\
  \num{9 }& \num{10893}   &\num{ 13} & \num{141609}  & ${5.324e\!-\!2}$&${ 1.467e\!-\!2 }$&${ 3.030e\!-\!2 }$&${ 4.498e\!-\!2 }$&${ 3.832e\!-\!2 }$&${ 8.331e\!-\!2 }$&${ 1.57 }$& \num{16 }\\
  \num{10} &\num{ 23386}  &\num{ 14} & \num{327404}  & ${3.854e\!-\!2}$&${ 5.454e\!-\!3 }$&${ 8.839e\!-\!3 }$&${ 1.429e\!-\!2 }$&${ 3.284e\!-\!2 }$&${ 4.713e\!-\!2 }$&${ 1.22 }$& \num{32 }\\
  \num{11} &\num{ 50202}  &\num{ 15} & \num{753030}  & ${3.165e\!-\!2}$&${ 1.990e\!-\!3 }$&${ 2.052e\!-\!3 }$&${ 4.042e\!-\!3 }$&${ 2.947e\!-\!2 }$&${ 3.351e\!-\!2 }$&${ 1.06 }$& \num{32 }\\
  \num{12} &\num{ 107461} &\num{ 16} & \num{1719376} & ${2.502e\!-\!2}$&${ 1.289e\!-\!3 }$&${ 1.360e\!-\!3 }$&${ 2.649e\!-\!3 }$&${ 2.765e\!-\!2 }$&${ 3.030e\!-\!2 }$&${ 1.21 }$& \num{32 }\\
  \num{13} &\num{ 230906} &\num{ 17} & \num{3925402} & ${2.305e\!-\!2}$&${ 7.541e\!-\!4 }$&${ 7.414e\!-\!4 }$&${ 1.495e\!-\!3 }$&${ 2.459e\!-\!2 }$&${ 2.608e\!-\!2 }$&${ 1.13 }$& \num{32 }\\
  \num{14} &\num{ 495538} &\num{18 }&\num{  8919684} & ${2.229e\!-\!2}$&${  4.669e\!-\!4}$&${  4.442e\!-\!4}$&${  9.111e\!-\!4}$&${  2.389e\!-\!2}$&${  2.480e\!-\!2}$&${  1.11}$& \num{32 }\\
  \num{15} &\num{1080281} &\num{19 }&\num{ 20525339} & ${2.131e\!-\!2}$&${  3.170e\!-\!4}$&${  3.125e\!-\!4}$&${  6.296e\!-\!4}$&${  2.269e\!-\!2}$&${  2.332e\!-\!2}$&${  1.09}$& \num{32 }\\
  \num{16} &\num{2333766} &\num{20 }&\num{ 46675320} & ${1.977e\!-\!2}$&${  2.067e\!-\!4}$&${  2.053e\!-\!4}$&${  4.121e\!-\!4}$&${  2.105e\!-\!2}$&${  2.147e\!-\!2}$&${  1.09}$& \num{64 }\\
  \num{17} &\num{5072224} &\num{22 }&\num{111588928} & ${1.495e\!-\!2}$&${  1.544e\!-\!4}$&${  1.560e\!-\!4}$&${  3.104e\!-\!4}$&${  1.789e\!-\!2}$&${  1.820e\!-\!2}$&${  1.21}$& \num{128}\\
  \bottomrule
\end{tabular}
\end{center}
\end{table}
%

\FloatBarrier%

\subsection{Nonstationary Hemker Problem}
\label{sec:Hemker}

\begin{figure}[hbt]
\centering
\begin{tikzpicture}[scale=0.66,thick]
\draw[black] (3,-3) -- ++(0,6);%
\draw[black] (-3, 3) -- ++ (6,-6);
\draw[black] (-3, -3) -- ++ (6,6);
\draw[black] (-3, 0) -- ++ (11,0);
\draw[black] (0, -3) -- ++ (0,6);
\draw[blue] (-3, -3) -- ++ (11,0);
\draw[blue] (8, -3) -- ++ (0,6);
\draw[blue] (-3, 3) -- ++ (11,0);
\draw[red,-latex] (2.5, 1.5) --node[above,red]{$\mathbf{b}$} ++ (2,0);
\draw[red,-latex] (2.5, -1.5) -- ++ (2,0);
\draw[green,fill=white] (0,0) circle (1);%
\draw[green] (-3,-3) -- node[left]{\color{green}$\Gamma_D$}++(0,6);%
\draw [magenta, fill] (-0.7071067811865475,0.7071067811865475) circle (3pt) node
[below right] {$\boldsymbol
  x_{e_2}$};
\draw [magenta, fill] (4,1) circle (3pt) node [right] {$\boldsymbol
  x_{e_1}$};
\end{tikzpicture}
\caption{\label{fig:hemker-geo}Geometry and coarse mesh of the domain for
the Hemker problem (left) and the best adaptive solution obtained in this
work (right). In the sketch of the geometry, green coloring corresponds to
Dirichlet BCs, while blue indicates homogeneous Neumann BCs. On the
circle, inhomogeneous Dirichlet BCSs are prescribed. The left boundary is
associated with homogeneous Dirichlet BCs.}
\end{figure}
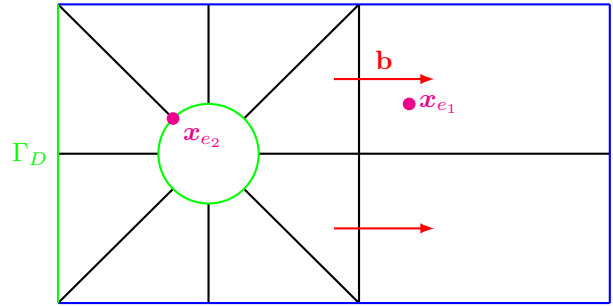

\begin{figure}[hbt]
  \centering
  \includegraphics[width=\linewidth]{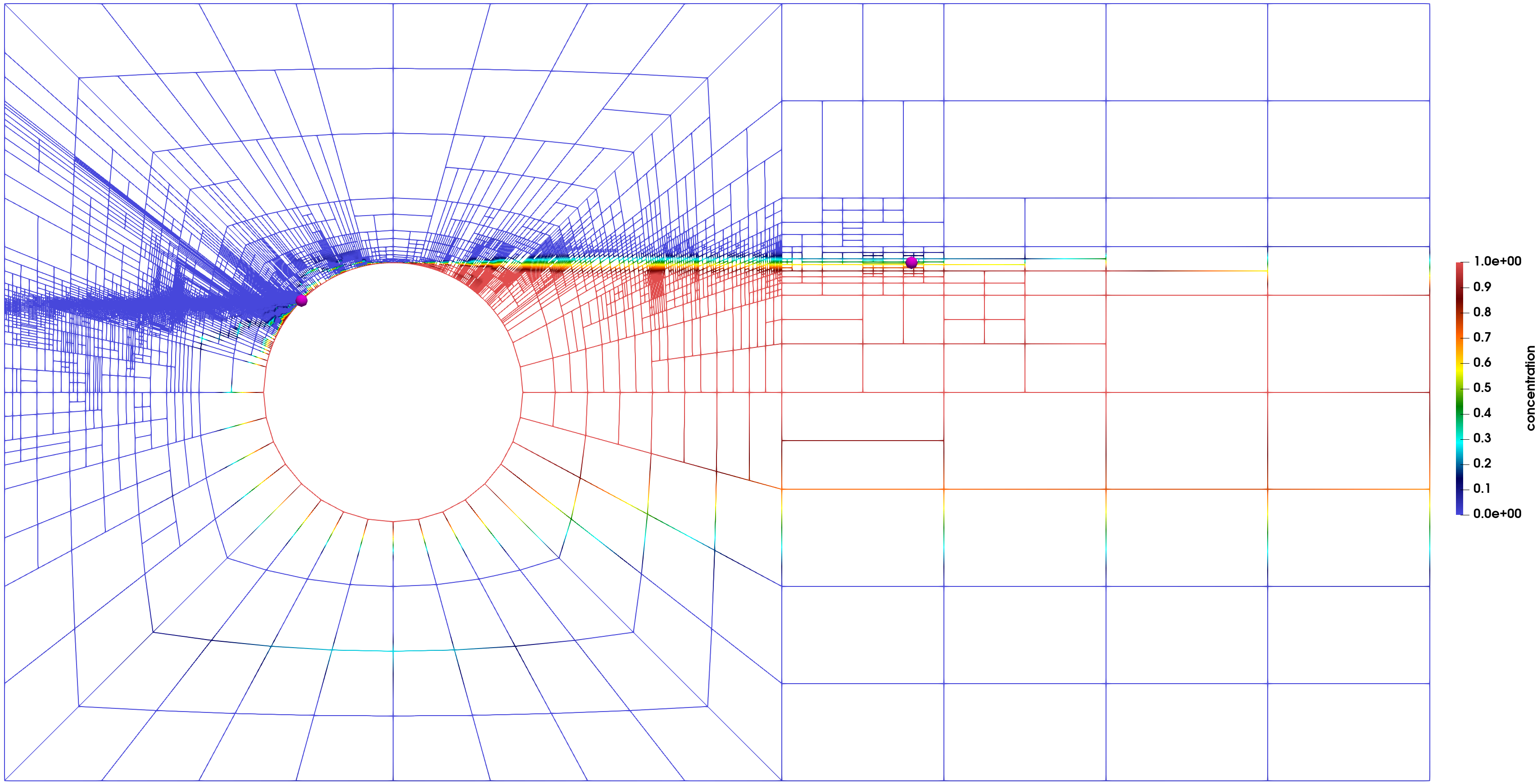}
  \caption{\label{fig:hemker-snapshot}The solution in the final DWR
    loop for Example~\ref{sec:Hemker}. One can clearly see the adaptive mesh refinement at
    $\boldsymbol{x}_{e_1}$ (right) and $\boldsymbol{x}_{e_2}$ (left) and upstream of the points due to the multi-goal~\eqref{eq:point-goal-1}.}
\end{figure}

\begin{figure}[hbt]
  \centering
  \includegraphics[width=.8\linewidth]{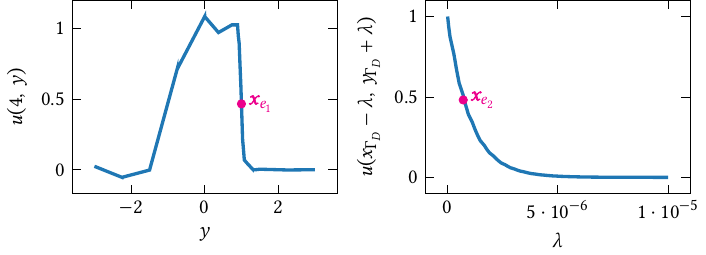}
  \caption{\label{fig:cut-lines}Cut lines of the solution to the nonstationary Hemker problem at
    final time. On the left, we
    plot a cut through the interior layers. On the right, we plot along the
    boundary normal at $(-2^{-1/2},\,2^{-1/2})$.
    The point of interest is marked in both plots, around which the smoothness
    of the solution can be observed. On the right, we see that the mesh is
    highly refined locally, allowing the notoriously difficult boundary layer to
    be resolved without over- or undershoots. Indeed, the solution near the goal
    point $\boldsymbol x_{e_2}$ is strictly positive between 0 and 1. }
\end{figure}

As a second example, we consider a modified version of the classical Hemker
problem introduced in~\cite{hemker_singularly_1996}, which is a more
sophisticated benchmark compared to the first example~\ref{sec:MovingHump}.
The problem is defined on the computational space-time domain given by
$\mathcal{Q}=\Omega\times I = \left( (-3,\,8) \times (-3,\,3) \right) \setminus \{ (x,\,y) \vert x^2 + y^2 \leq 1 \}\times (0,9]$, illustrated with the coarse spatial mesh given in in Figure~\ref{fig:hemker-geo}.
The boundaries colored green in Figure~\ref{fig:hemker-geo} correspond to Dirichlet boundary conditions.
Specifically, the boundary at \(x = -3\) is subjected to homogeneous Dirichlet
conditions $u=0$. The boundaries colored blue indicate homogeneous Neumann
boundary conditions $\partial_{\mathbf{n}} u=0$. On the circular boundary,
inhomogeneous Dirichlet boundary conditions are imposed with \(u = 1\). The
refinement fraction is set to
\(\theta_{\text{space}}^{\text{ref}} = \frac{1}{3}\), and the coarsening
fraction is \(\theta_{\text{space}}^{\text{co}} = 0\).
Here, we choose $\varepsilon = 10^{-6}$, $\boldsymbol b=(1,\,0)^\top$, $\alpha = 0$, as well as $f = 0$.
The aim is to control the error within two spatial control points
\begin{equation}
  \label{eq:hemker-control-points}
  \boldsymbol x_{e_1}=(4,\,1)\quad\text{and}\quad \boldsymbol x_{e_2}=(x_{e_2},\,y_{e_2})=(-2^{-1/2}-10^{-6},\,2^{-1/2}+10^{-6})\,,
\end{equation}
respectively,
which are located within the upper interior layer and the boundary layer at the
obstacle, respectively.
Corresponding to Sec.~\ref{sec:4:2:multi-goal}, we use the following multi-goal
\begin{equation}\label{eq:point-goal-1}
\boldsymbol J(u)=\big(J_1(u),\,J_2(u)\big)\,,\textnormal{ with } \quad
J_1(u)=\displaystyle\int_I u(\boldsymbol x_{e_1},\,t)\,\mathrm{d}t\,, \quad  J_2(u)=\displaystyle\int_I u(\boldsymbol x_{e_2},\,t)\,\mathrm{d}t\,.
\end{equation}
Here, the weights are set to one ($w_i=1\,,i\in\{1,2\}$) such that
\[J_w(u)=\displaystyle\sum_{i=1}^2J_i(u)\,.\]
Furthermore, we regularize
$J_i$, $i\in\{1,2\}$ analogously to Example~\ref{sec:MovingHump}. Moreover, we denote the projection of $\boldsymbol x_{e_2}$ onto the circular part of $\Gamma_D$ as $\boldsymbol x_{\Gamma_D}=(x_{\Gamma_D}\,y_{\Gamma_D})=(-2^{-1/2},\,2^{-1/2})$.

This test case is solved using a cG($1$)-dG($0$) discretization for the primal problem and a cG($2$)-dG($0$) discretization for the adjoint problem.
Figure~\ref{fig:hemker-snapshot} illustrates the solution and demonstrates
the effectiveness of anisotropic multi-goal oriented adaptive mesh refinement.
We observe anisotropic mesh refinement along the interior and boundary layer parallel to the
convection field to control the point-values in $\boldsymbol x_{e_1}$ and
$\boldsymbol x_{e_1}$, respectively. Cells without strong impact on the
multi-goal~\eqref{eq:point-goal-1} remain coarse. This holds in particular for
cells downstream of the goal points and in the lower half of the domain, where
we do not place any goal points. In the boundary layer, we observe anisotropic
mesh refinement mostly in $y$-direction at the top of the circular obstacle,
where the convection field is tangential to the boundary. We note that some of
the anisotropic refinement may result from anisotropic mesh smoothing. Overall
we obtain a very economic mesh with respect to the multi goal functional.

Fig.~\ref{fig:cut-lines} further underlines the effectiveness of our approach. We plot the solutions along cut lines through the two goal
points and interior and boundary layers. We observe that the upper interior
layer, where the goal point $\boldsymbol x_{e_1}$ lies, is well resolved while
the lower interior layer still exhibits oscillations. A similar behavor is
observed in the boundary layer close to $\boldsymbol x_{e_2}$. This is expected
and desired, since the chosen multi goal functional~\eqref{eq:point-goal-1} aims
to control the error in the points $\boldsymbol x_{e_1}$ and
$\boldsymbol x_{e_2}$ only.

Resolving the boundary layer at $\boldsymbol x_{e_2}$ requires a locally highly
refined mesh. Our results could provide a starting point to obtain reference values for the
thickness of the boundary layer in the Hemker example for $\varepsilon
=10^{-6}$. To the best of our knowledge, such results have not been documented
in the available literature on this particular setting.
We define the width of the boundary layer,
\begin{equation}\label{eq:ylayer}
  y_{\text{layer}}=\lVert \boldsymbol x_1-\boldsymbol x_0 \rVert\,,
\end{equation}
as the distance between the points $\boldsymbol x_1$ and $\boldsymbol x_0$ along
the normal of $\Gamma_D$ at $\boldsymbol x_{\Gamma_D}$. Here,
$\boldsymbol x_1\coloneq\boldsymbol x_{\Gamma_D}+
\lambda_1(-\frac{1}{\sqrt{2}},\frac{1}{\sqrt{2}})$ is defined as the point where
$u(\boldsymbol x_1) = 0.9$. Similarly,
$\boldsymbol x_0\coloneq\boldsymbol
x_{\Gamma_D}+\lambda_0(-\frac{1}{\sqrt{2}},\frac{1}{\sqrt{2}})$ is defined as
the point where $u(\boldsymbol x_0) = 0.1$. In the $22$nd DWR loop we obtain
$y_{\text{layer}}=\num{3.12e-06}$, while in the $23$rd DWR loop we obtain
$y_{\text{layer}}=\num{3.20e-06}$. Therefore, we obtain an almost converged
solution, which could be a starting point for further investigated. We note that
the problem remains solvable, but its size is already substantial for a
two-dimensional problem.
\begin{table}[htb]\setlength{\tabcolsep}{8pt}
  \caption{\label{tab:dwr-hemker-estimators}Anisotropic adaptive refinement including directional error indicators for goal quantity~\eqref{eq:point-goal-1} and maximum aspect ratios for Example~\ref{sec:Hemker}, $\varepsilon=10^{-6}$, $\delta_0=0.1$.}
  \begin{center}\footnotesize
  \begin{tabular}{rrrrrrrrrr}
    \toprule
    $\ell$ & $N_{\text{space}}$
    & $N_{t}$ & $N_{\text{tot}}$
    & $\eta_{h,x}$
    & $\eta_{h,y}$  & $\eta_h$ & $\eta_\tau$
    & $\eta_{\tau h}$
    & $\operatorname{ar}_{\max}$\\
    \midrule
      $ 1$& \num{    196}& $36$& \num{     7056}&  ${-2.16e\!+\!0}$& ${-2.98e\!+\!0}$& ${-5.14e\!+\!0}$& ${-2.25e\!-\!02}$& ${-5.16e\!+\!0}$&  ${   3.7}$\\
      $ 2$& \num{    310}& $36$& \num{    11160}&  ${-2.43e\!+\!0}$& ${-2.74e\!+\!0}$& ${-5.18e\!+\!0}$& ${-2.23e\!-\!03}$& ${-5.18e\!+\!0}$&  ${   5.3}$\\
      $ 3$& \num{    493}& $36$& \num{    17748}&  ${-2.67e\!+\!0}$& ${-3.09e\!+\!0}$& ${-5.76e\!+\!0}$& ${ 1.11e\!-\!03}$& ${-5.76e\!+\!0}$&  ${   7.2}$\\
      $ 4$& \num{    794}& $36$& \num{    28584}&  ${-2.05e\!+\!0}$& ${-2.06e\!+\!0}$& ${-4.12e\!+\!0}$& ${ 2.06e\!-\!02}$& ${-4.10e\!+\!0}$&  ${  14.4}$\\
      $ 5$& \num{   1288}& $36$& \num{    46368}&  ${-1.56e\!+\!0}$& ${-2.57e\!+\!0}$& ${-4.13e\!+\!0}$& ${ 4.36e\!-\!03}$& ${-4.12e\!+\!0}$&  ${  28.8}$\\
      $ 6$& \num{   2065}& $36$& \num{    74340}&  ${-2.18e\!+\!0}$& ${-2.25e\!+\!0}$& ${-4.43e\!+\!0}$& ${ 1.91e\!-\!02}$& ${-4.41e\!+\!0}$&  ${  29.6}$\\
      $ 7$& \num{   3262}& $36$& \num{   117432}&  ${-1.34e\!+\!0}$& ${-2.63e\!+\!0}$& ${-3.98e\!+\!0}$& ${ 2.27e\!-\!03}$& ${-3.98e\!+\!0}$&  ${  29.6}$\\
      $ 8$& \num{   5246}& $36$& \num{   188856}&  ${-1.66e\!+\!0}$& ${-2.73e\!+\!0}$& ${-4.40e\!+\!0}$& ${-4.84e\!-\!05}$& ${-4.40e\!+\!0}$&  ${  40.7}$\\
      $ 9$& \num{   8095}& $36$& \num{   291420}&  ${-1.31e\!+\!0}$& ${-2.69e\!+\!0}$& ${-4.01e\!+\!0}$& ${ 4.05e\!-\!05}$& ${-4.01e\!+\!0}$&  ${  49.7}$\\
      $10$& \num{  12663}& $36$& \num{   455868}&  ${-1.34e\!+\!0}$& ${-2.72e\!+\!0}$& ${-4.06e\!+\!0}$& ${-1.33e\!-\!05}$& ${-4.06e\!+\!0}$&  ${  93.3}$\\
      $11$& \num{  19967}& $36$& \num{   718812}&  ${-1.61e\!+\!0}$& ${-2.53e\!+\!0}$& ${-4.14e\!+\!0}$& ${ 9.37e\!-\!06}$& ${-4.14e\!+\!0}$&  ${ 186.7}$\\
      $12$& \num{  31629}& $36$& \num{  1138644}&  ${-1.05e\!+\!0}$& ${-1.06e\!+\!0}$& ${-2.11e\!+\!0}$& ${ 5.52e\!-\!06}$& ${-2.11e\!+\!0}$&  ${ 186.7}$\\
      $13$& \num{  50214}& $36$& \num{  1807704}&  ${-4.95e\!-\!1}$& ${-5.01e\!-\!1}$& ${-9.97e\!-\!1}$& ${ 1.49e\!-\!06}$& ${-9.97e\!-\!1}$&  ${ 194.5}$\\
      $14$& \num{  80057}& $36$& \num{  2882052}&  ${-2.46e\!-\!1}$& ${-2.47e\!-\!1}$& ${-4.94e\!-\!1}$& ${ 2.96e\!-\!07}$& ${-4.94e\!-\!1}$&  ${ 353.2}$\\
      $15$& \num{ 127751}& $36$& \num{  4599036}&  ${-1.23e\!-\!1}$& ${-1.23e\!-\!1}$& ${-2.46e\!-\!1}$& ${ 4.69e\!-\!08}$& ${-2.46e\!-\!1}$&  ${ 374.8}$\\
      $16$& \num{ 203497}& $36$& \num{  7325892}&  ${-6.16e\!-\!2}$& ${-6.20e\!-\!2}$& ${-1.23e\!-\!1}$& ${ 4.51e\!-\!09}$& ${-1.23e\!-\!1}$&  ${ 375.1}$\\
      $17$& \num{ 324594}& $36$& \num{ 11685384}&  ${-3.08e\!-\!2}$& ${-3.14e\!-\!2}$& ${-6.23e\!-\!2}$& ${-2.60e\!-\!09}$& ${-6.23e\!-\!2}$&  ${ 434.5}$\\
      $18$& \num{ 517116}& $36$& \num{ 18616176}&  ${-1.54e\!-\!2}$& ${-1.63e\!-\!2}$& ${-3.17e\!-\!2}$& ${-2.02e\!-\!09}$& ${-3.17e\!-\!2}$&  ${ 745.7}$\\
      $19$& \num{ 823827}& $36$& \num{ 29657772}&  ${-7.70e\!-\!3}$& ${-8.68e\!-\!3}$& ${-1.63e\!-\!2}$& ${-1.33e\!-\!09}$& ${-1.63e\!-\!2}$&  ${ 764.6}$\\
      $20$& \num{1312642}& $36$& \num{ 47255112}&  ${-3.85e\!-\!3}$& ${-4.69e\!-\!3}$& ${-8.54e\!-\!3}$& ${-7.47e\!-\!10}$& ${-8.54e\!-\!3}$&  ${1529.1}$\\
      $21$& \num{2070136}& $36$& \num{ 74524896}&  ${-1.92e\!-\!3}$& ${-2.52e\!-\!3}$& ${-4.45e\!-\!3}$& ${-3.97e\!-\!10}$& ${-4.45e\!-\!3}$&  ${1529.1}$\\
      $22$& \num{3194060}& $36$& \num{114986160}&  ${-9.63e\!-\!4}$& ${-1.33e\!-\!3}$& ${-2.30e\!-\!3}$& ${-2.06e\!-\!10}$& ${-2.30e\!-\!3}$&  ${1529.1}$\\
      $23$& \num{4942506}& $36$& \num{177930216}&  ${-4.81e\!-\!4}$& ${-6.92e\!-\!4}$& ${-1.17e\!-\!3}$& ${-1.05e\!-\!10}$& ${-1.17e\!-\!3}$&  ${1529.1}$\\
    \bottomrule
  \end{tabular}
\end{center}
\end{table}

In Tab.~\ref{tab:dwr-hemker-estimators} we show the evolution of the spatial
and temporal error estimators over the DWR loops. We clearly observe the
convergence of the error estimators. The maximum aspect ratio of the spatial
mesh increases from 3.7 to 1529.1, which underscores the effectiveness of
anisotropic mesh refinement.

\FloatBarrier%
\section{Conclusion}
In this work, we presented an anisotropic multi-goal error control strategy for
time-dependent convection-diffusion-reaction equations based on the Dual
Weighted Residual (DWR) method. By employing anisotropic interpolation and
restriction operators, the proposed method delivers space- and time-resolved
error indicators, which additionally decouple the spatial contributions
directionally. These directional indicators quantify anisotropy of the solution
with respect to multiple goal functionals and enable the construction of highly
anisotropic meshes that accurately resolve sharp layers.

To suppress spurious oscillations in convection-dominated regimes, SUPG
stabilization is incorporated in both primal and adjoint problems. The resulting
adaptive algorithm reliably generates high-aspect-ratio elements, in regions
where they are most effective. This enables efficient resolution of interior and
boundary layers without producing numerical artifacts.

Numerical experiments demonstrate the accuracy and robustness of the method
across a range of multi-goal functionals. An analytical test case confirms the
effectivity and reliability of the proposed error estimators, with observed
effectivity indices close to one. The challenging Hemker problem is successfully
addressed, and the adaptive strategy enables efficient resolution of the
boundary layer at a diffusion coefficient of $\varepsilon=10^{-6}$. Initial
steps are taken toward generating reliable reference values for the boundary
layer thickness in this convection-dominated regime.

Due to its generic structure, the proposed method can be straightforwardly
applied to a wide range of convection-dominated problems with multiple goals.

\begin{acknowledgement}
All authors acknowledge the funding of DAAD-project 57729992,
"Goal-oriented AnIsotropic Space-Time Mesh Adaption (AIMASIM)" in the funding program "Programm des projektbezogenen Personenaustauschs Griechenland ab 2024".
Bernhard Endtmayer and Thomas Wick additionally acknowledge the support by the Cluster of Excellence PhoenixD (EXC 2122, Project ID 390833453).
Bernhard Endtmayer was funded by an Humboldt Postdoctoral Fellowship at the beginning of the work.
Bernhard Endtmayer thanks the research group of Markus Bause for financing the research visit at the Helmut Schmidt University Hamburg in February 16--28, 2025.
Computational resources (HPC cluster HSUper) have been provided by the project hpc.bw, funded by dtec.bw - Digitalization and Technology Research Center of the Bundeswehr. dtec.bw is funded by the European Union - NextGenerationEU.
\end{acknowledgement}

\vspace{\baselineskip}

\end{document}